\documentclass[
11pt]{amsart}

\usepackage[latin1]{inputenc}
\usepackage{subfigure}

\usepackage[round]{natbib}

\usepackage{amsmath}
\usepackage{amsfonts}
\usepackage{amssymb}
\usepackage{amsthm}
\usepackage{color}
\usepackage{enumerate}
\usepackage{amscd}
\usepackage{epsfig}
\usepackage{multicol}
\usepackage{xypic}

\newcommand{\comment}[1]{}

\theoremstyle{remark}
\newtheorem{Rem}{Remark}
\newtheorem{example}{Example}

\newtheorem{Def}{Definition}[section]

\newtheorem{Th}{Theorem}[section]
\newtheorem{Defth}[Th]{Definition-Theorem}
\newtheorem{Prop}[Th]{Proposition}
\newtheorem{Lem}[Th]{Lemma}
\newtheorem{Corol}[Th]{Corollary}

\newcommand{\pos}{\mathcal{P}}
\newcommand{\ext}{\mathcal{L}}
\newcommand{\sumext}{\varphi}

\newcommand{\poset}{{\tt poset}}

\DeclareMathOperator{\extr}{ext}

\author[A. Boussicault and V. Féray]{Adrien Boussicault \and Valentin Féray}
\title[Rational identities and graph combinatorics]{Application of graph combinatorics to rational identities of type $A$}
\address{
Université Paris-Est,
Institut d'électronique et d'informatique Gaspard-Monge,
77454 Marne-la-Vallée Cedex 2\\
}
\email{feray@univ-mlv.fr, boussica@univ-mlv.fr} 

\input xy
\xyoption{all}

\begin{document}

\maketitle

\begin{abstract}

To a word $w$, we associate the rational function $\Psi_w = \prod (x_{w_i} - x_{w_{i+1}})^{-1}$. The main object, introduced by C. Greene to generalize identities linked to Murnaghan-Nakayama rule, is a sum of its images by certain permutations of the variables. The sets of permutations that we consider are the linear extensions of oriented graphs. We explain how to compute this rational function, using the combinatorics of the graph $G$. We also establish a link between an algebraic property of the rational function (the factorization of the numerator) and a combinatorial property of the graph (the existence of a disconnecting chain).\\

\end{abstract}

\section{Introduction}

A partially ordered set (poset) $\mathcal{P}$ is a finite set $V$ endowed with a partial order. By definition, a word $w$ containing exactly once each element of $V$ is called a linear extension if the order of its letters is compatible with $\mathcal{P}$ (if $a \leq_{\mathcal{P}} b$, then $a$ must be before $b$ in $w$). To a linear extension $w=v_1 v_2 \ldots v_n$, we associate a rational function:
$$\psi_w = \frac{1}{(x_{v_1}-x_{v_2})\cdot(x_{v_2}-x_{v_3})\ldots (x_{v_{n-1}}-x_{v_n})}.$$

We can now introduce the main object of the paper. If we denote by $\ext(\mathcal{P})$ the set of linear extensions of $\mathcal{P}$, then we define $\Psi_{\mathcal{P}}$ by:
$$\Psi_{\mathcal{P}} = \sum_{w \in \ext(\mathcal{P})} \psi_w.$$

\subsection{Background}\label{subsect:background}

The linear extensions of posets contain very interesting subsets of the symmetric group: for example, the linear extensions of the poset considered in the article \cite{BUTLER} are the permutations smaller than a permutation $\pi$ for the weak Bruhat order. In this case, our construction is close to that of Demazure characters \cite{DEMAZURE}. S. Butler and M. Bousquet-M\'elou characterize the permutations $\pi$ corresponding to acyclic posets, which are exactly the cases where the function we consider is the simplest.\\

Moreover, linear extensions are hidden in a recent formula for irreducible character values of the symmetric group: if we use the notations of \cite{FS}, the quantity $N^\lambda(G)$ can be seen as a sum over the linear extensions of the bipartite graph $G$ (bipartite graphs are a particular case of oriented graphs). This explains the similarity of the combinatorics in article \cite{Fe} and in this one.\\

The function $\Psi_{\mathcal{P}}$ was considered by C. Greene \cite{GREENE}, who wanted to generalize a rational identity linked to Murnaghan-Nakayama rule for irreducible character values of the symmetric group. He has given in his article a closed formula for planar posets ($\mu_\pos$ is the M\"obius function of $\pos$):
$$
\Psi_{\mathcal{P}} =
\left\{
\begin{array}{cl}
0 & \text{if } \mathcal{P} \text{ is not connected,} \\
\prod\limits_{y,z \in \pos} (x_y-x_z)^{\mu_\pos(y,z)} & \text{if } \mathcal{P} \text{ is connected,}
\end{array}
\right.
$$
However, there is no such formula for general posets, only the denominator of the reduced form of $\Psi_{\mathcal{P}}$ is known \cite{Bo}. In this article, the first author has investigated the effects of elementary transformations of the Hasse diagram of a poset on the numerator of the associated rational function. He has also noticed, that in some case, the numerator is a specialization of a Schur function \cite[paragraph 4.2]{Bo} (we can also find multiSchur functions or Schubert polynomials).\\

In this paper, we obtain some new results on this numerator, thanks to a simple local transformation in the graph algebra, preserving linear extensions.

\subsection{Main results}
\subsubsection{An inductive algorithm}
The first main result of this paper is an induction relation on linear extensions (Theorem \ref{th:boucle}). When one applies $\Psi$ on it, it gives an efficient algorithm to compute the numerator of the reduced fraction of $\Psi_\mathcal{P}$ (the denominator is already known).

\subsubsection{A combinatorial formula}
If we iterate our first main result in a clever way, we can describe combinatorially the final result. The consequence is our second main result: if we give to the graph of a poset $\mathcal{P}$ a rooted map structure, we have a combinatorial non-inductive formula for the numerator of $\Psi_{\mathcal{P}}$ (Theorem \ref{th:N_combi}).

\subsubsection{A condition for $\Psi_\mathcal{P}$ to factorize}
Greene's formula for the function associated to a planar poset is a quotient of products of polynomials of degree $1$. In the non-planar case, the denominator is still a product of degree 1 terms, but not the numerator. So we may wonder when the numerator $N(\mathcal{P})$ can be factorized.\\
Our third main result is a partial answer (a sufficient but not necessary condition) to this question: the numerator $N(\mathcal{P})$ factorizes if there is a chain disconnecting the Hasse diagram of $\mathcal{P}$ (see Theorem \ref{factorization_chaine} for a precise statement). An example is drawn on figure \ref{fig:ex_factorization} (the disconnecting chain is $(2,5)$). Note that we use here and in the whole paper a unusual convention: we draw the posets from left (minimal elements) to right (maximal elements).

\begin{figure}[h]
\begin{center}
$
N \left(
	\begin{array}{c}
	\mbox{\epsfig{file=figure-1.ps,width=24 pt}}
	\end{array} \right)
=
N \left(
	\begin{array}{c}  
	\mbox{\epsfig{file=figure-2.ps,width=24 pt}}
	\end{array} \right)
.
N \left(
	\begin{array}{c}
	\mbox{\epsfig{file=figure-3.ps,width=24 pt}}
	\end{array} \right)
$
\caption{Example of chain factorization\label{fig:ex_factorization}}
\end{center}
\end{figure}

\subsection{Open problems}

\subsubsection{Around the map structure} \label{subsubsect:openpb_mapstructure}
Theorem \ref{th:N_combi} is a cominatorial formula for the numerator of $\Psi_\mathcal{P}$ involving a map structure on the corresponding graph. Can we find a formula, which does not depend any additional structure on the graph?

Furthermore if we use ordered-embeddings of graphs in $\mathbb{R} \times \mathbb{R}$ (see definition \ref{def:planar}), the map structure is not independant from the poset structure. Is there a way to use this link?

\subsubsection{Necessary condition for factorization}\label{subsubsect:openpb_facto}
The conclusion of the factorization Theorem \ref{factorization_chaine} is sometimes true, even when the separating path is not a chain: see for example Figure \ref{fig:other:factorization} (the path $(5,6,3)$ disconnects the Hasse diagram, but is not a chain).\\
This equality, and many more, can be easily proved using the same method as Theorem \ref{factorization_chaine}. Can we give a necessary (and sufficient) condition for the numerator of a poset to factorize into a product of numerators of subposets? Are all factorizations of this kind?

\begin{figure}
\begin{center}
$
N \left(
	\begin{array}{c}
	\mbox{\epsfig{file=figure-4.ps,width=60 pt}}
	\end{array} \right)
=
N \left(
	\begin{array}{c}
	\mbox{\epsfig{file=figure-5.ps,width=48 pt}}
	\end{array} \right)
.
N \left(
	\begin{array}{c}  
	\mbox{\epsfig{file=figure-6.ps,width=48 pt}}
	\end{array} \right)
$
\caption{An example of factorization, not contained in Theorem \ref{factorization_chaine}.  \label{fig:other:factorization}}
\end{center}
\end{figure}

\subsubsection{Characterisation of the numerator}\label{subsubsect:openpb_annulation}
Let us consider a bipartite poset $\mathcal{P}$ (which has only minimal and maximal elements, respectively $a_1,\ldots,a_l$ and $b_1,\ldots,b_r$). The numerator $N(\mathcal{P})$ of $\Psi_{\mathcal{P}}$ is a polynomial in $b_1,\ldots,b_r$ which degree in each variable can be easily bounded \cite[Proposition 3.1]{Bo}. Moreover, we know, by Corollary \ref{corol:annulation}, that $N(\mathcal{P})=0$ on some affine subspaces of the space of variables. Unfortunately, these vanishing relations and its degree do not characterize $N(\mathcal{P})$ up to a multiplicative factor. Is there a bigger family of vanishing relations, linked to the combinatorics of the Hasse diagram of the poset, which characterizes $N(\mathcal{P})$?\\
This question comes from the following observation: for some particular posets, the numerator is a Schubert polynomial and Schubert polynomials are known to be easily defined by vanishing conditions \cite{LASCOUX2008}.

\subsection{Outline of the paper}
In section \ref{sect:graphs}, we present some basic definitions on graphs and posets.\\

In section \ref{sect3}, we introduce our main object and its basic properties.\\

In section \ref{sect_oper}, we state our first main result: an inductive relation for linear extensions. The next section (\ref{sect_exemples}) is devoted to some explicit computations using this result.\\

Section \ref{sect:combi} gives a combinatorial description of the result of the iteration of our inductive relation: we derive from it our second main result, a combinatorial formula for the numerator of $\Psi_{\mathcal{P}}$.\\

The last Section (\ref{sectchainfact}) is devoted to our third main result: a sufficient condition of factorization.

\section{Graphs and posets} \label{sect:graphs}
Oriented graphs are a natural way to encode information of posets. To avoid confusions, we recall all necessary definitions in paragraph \ref{graphs_def}. The definition of linear extensions can be easily formulated directly in terms of graphs (paragraph \ref{subsect:posetgraphs}).\\
We will also define some elementary removal operations on graphs (paragraph \ref{operations_graph}), which will be used in the next section. Due to transitivity relations, it is not equivalent to perform these operations on the Hasse diagram or on the complete graph of a poset, that's why we prefer to formulate everything in terms of graphs.

\subsection{Definitions and notations on graphs}\label{graphs_def}

In this paper, we deal with finite \emph{directed graphs}. So we will use the following definition of a graph $G$:
\begin{itemize}
\item A finite set of vertices $V_G$.
\item A set of edges $E_G$ defined by $E_G \subset V_G \times V_G$.
\end{itemize}
If $e \in E_G$, we will note by $\alpha(e) \in V_G$ the first component of $e$ (called \emph{origin} of $e$) and $\omega(e)\in V_G$ its second component (called \emph{end} of $e$). This means that each edge has an orientation.\\
Let $e = (v_1,v_2)$ be an element of $V_G \times V_G$. Then we denote by $\overline{e}$ the pair $(v_2,v_1)$.

With this definition of graphs, we have four definitions of injective walks on the graph.

$$\begin{array}{c|c|c}
  & \begin{array}{c} \text{can not go } \text{backwards} \end{array} & \begin{array}{c}\text{can go } \text{backwards}\end{array} \\
\hline
\text{closed} & \text{circuit} & \text{cycle}\\
\hline
\text{non-closed} & \text{chain} & \text{path}
\end{array}$$

More precisely,
\begin{Def} Let $G$ be a graph and $E$ its set of edges.
\begin{description}
 \item[chain] A chain is a sequence of edges $c=(e_1, \hdots, e_k)$ of $G$ such that $\omega(e_1)=\alpha(e_2)$, $\omega(e_2)=\alpha(e_3)$, $\hdots$ and $\omega(e_{k-1})=\alpha(e_k)$.
 \item[circuit] A circuit is a chain $(e_1,\hdots ,e_k)$ of $G$ such that $\omega(e_k)=\alpha(e_1)$.
 \item[path] A path is a sequence $(e_1,\ldots,e_h)$ of elements of $E \cup \overline{E}$ such that $\omega(e_1)=\alpha(e_2)$, $\omega(e_2)=\alpha(e_3)$, $\hdots$ and $\omega(e_{k-1})=\alpha(e_k)$.
 \item[cycle] A cycle $C$ is a path with the additional property that $\omega(e_k)=\alpha(e_1)$. If $C$ is a cycle, then we denote by $HE(C)$ the set $C \cap E$.
\end{description}
In all these definitions, we add the condition that all edges and vertices are different (except of course, the equalities in the definition).
\end{Def}

\begin{Rem}
The difference between a cycle and a circuit (respectively a path and a chain) is that, in a cycle (respectively in a path), an edge can appear in both directions (not only in the direction given by the graph structure). The edges, which appear in a cycle $C$ with the same orientation than their orientation in the graph, are exactly the elements of $HE(C)$.
\end{Rem}

To make the figures easier to read, $\alpha(e)$ is always the left-most extremity of $e$ and $\omega(e)$ its right-most one.
Such drawing construction is not possible if the graph contains a circuit.
But its case will not be very interesting for our purpose.

\begin{example}
 An example of graph is drawn on figure \ref{fig:ex_graph}. In the left-hand side, the non-dotted edges form a chain $c$, whereas, in the right-hand side, they form a cycle $C$, such that $HE(C)$ contains 3 edges: $(1,6), (6,8)$ and $(5,7)$.
\begin{figure}
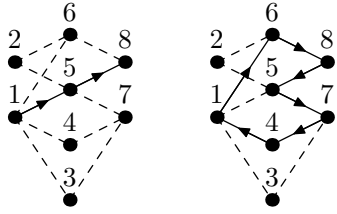

$$\epsfig{file=figure-7.ps,width=48 pt} \hspace{1cm} \epsfig{file=figure-8.ps,width=48 pt}
 $$
\caption{Example of a chain and a cycle $C$ (we recall that orientations are from left to right).}
\label{fig:ex_graph}
\end{figure}
\end{example}
\comment{
Let $e = (v_1,v_2)$ be an edge and $C=(e_1,e_2, \dots, e_n)$ be a cycle. We denote by $\overline{e}$ and $\overline{C}$ the edge $\overline{e}=(v_2,v_1)$ and the cycle $\overline{C} = \overline{e_n}, \dots, \overline{e_1}$ with $HE(\overline{C})= L \setminus\overline{HE(C)}:= \{ \overline{e} | e \in HE(C) \}$.
}

The {\it cyclomatic number} of a graph $G$ is $|E_G| - |V_G| + c_G$, where $c_G$ is the number of connected components of $G$.
A graph contains a cycle if and only if its cyclomatic number is not $0$ (see \cite{GRAPH}).
If it is not the case, the graph is called {\it forest}.
A connected forest is, by definition, a tree.
Beware that, in this context, there are no rules for the orientation of the edges of a tree (often, in the literature, an oriented tree is a tree which edges are oriented from the \emph{root} to the \emph{leaves}, but we do not consider such objects here).

\subsection{Posets, graphs, Hasse diagrams and linear extensions}\label{subsect:posetgraphs}

In this paragraph, we recall the link between graphs and posets.\\

Given a graph $G$, we can consider the binary relation on the set $V_G$ of vertices of $G$:
$$x \le y \stackrel{\text{\tiny def}}{\Longleftrightarrow} \left( x=y \text{ or } \exists \ e \in E_G \text{ such that } \left\{ \begin{array}{c} \alpha(e)=x \\ \omega(e)=y \end{array} \right. \right)$$

This binary relation can be completed by transitivity. If the graph has no circuit, the resulting relation $\le$ is antisymmetric and, hence, endows the set $V_G$ with a poset structure, which will be denoted $\poset(G)$.\\

\begin{figure}
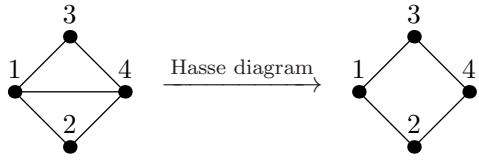

\begin{center}
\begin{center}
$
\begin{CD}
\begin{array}{c}
\mbox{\epsfig{file=figure-9.ps,width=48 pt}}
\end{array}
@>\text{Hasse diagram}>>
\begin{array}{c} 
\mbox{\epsfig{file=figure-10.ps,width=48 pt}}
\end{array}\\
\end{CD}
$
\end{center}
\caption{Example of a poset and his Hasse diagram.}
\label{fig_hasse_diagram}
\end{center}
\end{figure}

The application $\poset$ is not injective. Among the pre-images of a given poset $\mathcal{P}$, there is a minimum one (for the inclusion of edge set), which is called Hasse diagram of $\mathcal{P}$ (see figure \ref{fig_hasse_diagram} for an example).\\

The definition of linear extensions given in the introduction can be formulated in terms of graphs:
\begin{Def}
A linear extension of a graph $G$ is a total order $\leq_w$ on the set of vertices $V$ such that,
for each edge $e$ of $G$, one has ${\alpha(e)} \leq_w {\omega(e)}$.\\

The set of linear extensions of $G$ is denoted $\ext(G)$. Let us also define the formal sum $\sumext(G)=\sum\limits_{w \in \ext(G)} w$.
\end{Def}

We will often see a total order $\leq_w$ defined by $v_{i_1} \leq_w v_{i_2} \leq_w \ldots \leq_w v_{i_n}$ as a word $w = v_{i_1} v_{i_2} \ldots v_{i_n}$.

For example, the linear extensions of the poset drawn in the figure \ref{fig_hasse_diagram} are $1234$ and $1324$.

\begin{Rem}\label{rem_circ}
If $G$ contains a circuit, then it has no linear extensions. Else, its linear extensions are the linear extensions of $\poset(G)$. Thus considering graphs instead of posets does not give more general results.\end{Rem}

The following lemma comes straight forward from the definition:
\begin{Lem}\label{lem_ext}
Let $G$ and $G'$ be two graphs with the same set of vertices. Then one has:
$$E(G) \subset E(G') \text{ and } w \in \ext(G') \Longrightarrow w \in \ext(G);$$
$$w \in \ext(G) \text{ and } w \in \ext(G') \Longleftrightarrow w \in \ext(G \vee G'),$$
where $G \vee G'$ is defined by $\left\{\begin{array}{l} V(G \vee G') = V(G) = V(G');\\E(G \vee G')=E(G) \cup E(G'). \end{array}\right.$
\end{Lem}

\subsection{Elementary operations on graphs}\label{operations_graph}
The main tool of this paper consists in removing some edges of a graph $G$.
\begin{Def} \label{def:remove_edges}
Let $G$ be a graph and $E'$ a subset of its set of edges $E_G$. We will denote by $G \backslash E'$ the graph $G'$ with
\begin{itemize}
 \item the same set of vertices as $G$ ;
 \item the set of edges $E_{G'}$ defined by $E_{G'}:=E_G \backslash E'$.
\end{itemize}
\end{Def}

\begin{Def}
 If $G$ is a graph and $V'$ a subset of its set of vertices $V$, $V'$ has an induced graph structure: its edges are exactly the edges of $G$, which have both their extremities in $V'$.\\

\begin{figure}
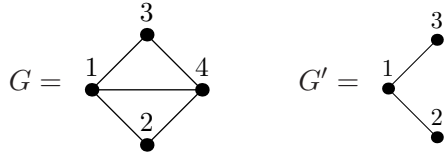

\begin{center}
$
G=
\begin{array}{c}
\mbox{\epsfig{file=figure-9.ps,width=48 pt}}
\end{array}
\hspace{1cm}
G'=
\begin{array}{c} 
\mbox{\epsfig{file=figure-11.ps,width=24 pt}}
\end{array}
$
\end{center}
\caption{$G'$ is the induced graph of $G$ by $\{ 1,2,3 \}$.}
\label{fig_induction}
\end{figure}

If $V \backslash V'=\{v_1,\ldots,v_l\}$, the graph induced by $V'$ will be denoted by $G \backslash \{v_1,\ldots, v_l\}$.
The symbol is the same than in definition \ref{def:remove_edges}, but it should not be confusing.
\end{Def}

\begin{Def}[Contraction] \label{def:contraction}
We denote by $G/e$ the graph (here, the set of edges can be a multiset) obtained by contracting the edge $e$ (\textit{i.e.} in $G/e$, there is only one vertex $v$ instead of $v_1$ and $v_2$, the edges of $G$ different from $e$ are edges of $G/e$: if their origin and/or end in $G$ is $v_1$ or $v_2$, it is $v$ in $G/e$).
\end{Def}

Then, if $\alpha(e) \neq \omega(e)$, $G/e$ is a graph with the same number of connected components and the same cyclomatic number as $G$.\\

\begin{figure}
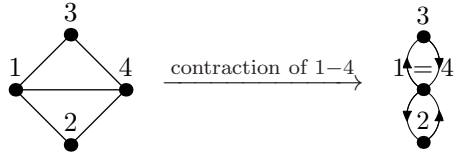

\begin{center}
$
\begin{CD}
\begin{array}{c}
\mbox{\epsfig{file=figure-9.ps,width=48 pt}}
\end{array}
  @>\text{contraction of }1-4>> 
\begin{array}{c} 
\mbox{\epsfig{file=figure-12.ps,width=24 pt}}
\end{array}
\end{CD}
$
\end{center}
\caption{Example of contraction.}
\label{fig_contraction}
\end{figure}

These definitions are illustrated on figures \ref{fig_induction} and \ref{fig_contraction}.

\section{Rational functions on graphs}\label{sect3}

\subsection{Definition}
Given a graph $G$ with $n$ vertices $v_1,\ldots,v_n$, we are interested in the following rational function $\Psi(G)$ in the variables $(x_{v_i})_{i=1\ldots n}$:
$$\Psi(G)=\sum_{w \in \ext(G)}\frac{1}{(x_{w_1} -x_{w_2}) \ldots (x_{w_{n-1}} - x_{w_n})}.$$

We also consider the renormalization:
$$N(G):= \Psi(G) \cdot \prod_{e \in E_G} ( x_{\alpha(e)} - x_{\omega(e)} ).$$
In fact, we will see later that it is a polynomial. Moreover, if $G$ is the Hasse diagram of a poset, $\displaystyle \Psi(G)=\frac{N(G)}{\prod\limits_{e \in E_G} ( x_{\alpha(e)} - x_{\omega(e)} )}$ is a reduced fraction.

\subsection{Pruning invariance}
Thanks to the following lemma, it will be easy to compute $N$ on forests (note that these results have already been proved in \cite{Bo}, but the following demonstrations are simpler and make this article self-contained).
\begin{Lem}\label{lem:emond}
Let $G$ be a graph with a vertex $v$ of valence $1$ and $e$ the edge of extremity (origin or end) $v$. Then one has
$$N(G) = N \big( G \backslash \{v\} \big).$$
\end{Lem}

For example,
$$
N \left( 
\begin{array}{c} \mbox{\epsfig{file=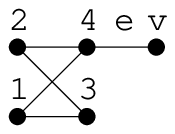,width=48 pt}} \end{array}
\right) =
N \left(
\begin{array}{c} \mbox{\epsfig{file=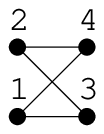,width=24 pt}} \end{array}
 \right) = x_1+x_2-x_3-x_4.
$$

\begin{proof}
One wants to prove that:
$$(x_{\alpha(e)}-x_{\omega(e)})\cdot \left(\sum_{w' \in \ext(G)} \psi_{w'}\right)=\sum_{w \in \ext(G \setminus \{v\})} \psi_w.$$
But one has a map $\text{er}_v : \ext(G) \rightarrow \ext(G \setminus \{v\})$ which sends a word $w'$ to the word $w$ obtained from $w'$ by erasing the letter $v$ (see figure \ref{fig_example_er}). So it is enough to prove that, for each $w \in \ext(G \setminus \{v\})$, one has :
$$(x_{\alpha(e)}-x_{\omega(e)})\cdot \left(\sum_{w' \in \text{er}_v^{-1}(w)} \psi_{w'}\right)=\psi_{w}.$$
Let us assume that $v$ is the end of $e$ and $w=w_1 \ldots w_{n-1} \in \ext(G \setminus \{v\})$. We denote by $k$ the index in $w$ of the origin of $e$. The set $\text{er}_v^{-1}(w)$ is: 
$$\big\{w_1 \ldots w_i v w_{i+1} \ldots w_{n-1}, i \geq k \big\}$$

\begin{figure}
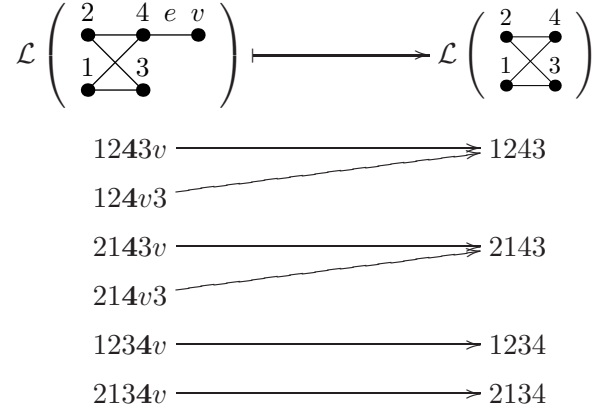

$$
\begin{array}{c}
\xymatrix@R=5pt@C=30pt{
\mathcal{L} \left(\begin{array}{c} \mbox{\epsfig{file=figure-13.ps,width=48 pt}} \end{array} \right) \ar@{|->}[rr]  & & \mathcal{L} \left( \begin{array}{c} \mbox{\epsfig{file=figure-14.ps,width=24 pt}} \end{array} \right) \\
12\mathbf{4}3v  \ar@{->}[rr] & &1243\\
12\mathbf{4}v3  \ar@{->}[urr] & &    \\
21\mathbf{4}3v  \ar@{->}[rr] & &2143\\
21\mathbf{4}v3  \ar@{->}[urr] & &    \\
123\mathbf{4}v  \ar@{->}[rr] & &1234\\
213\mathbf{4}v  \ar@{->}[rr] & &2134
}
\end{array}
$$
\caption{Example of the map $er_v$}
\label{fig_example_er}
\end{figure}

So, one has:
\begin{eqnarray*}
\sum_{w' \in \text{er}_v^{-1}(w)} \psi_{w'}  &=& \sum_{i=k}^{n-1} \frac{1}{\left[\begin{array}{l} (w_1 - w_2) \ldots (w_{i-1} - w_i) (w_i - v) \\ \qquad \cdot (v-w_{i+1}) (w_{i+1} - w_{i+2}) \ldots (w_{n-2} - w_{n-1}) \end{array}\right]} \\
&=& \frac{1}{
\left[
\begin{array}{l}
(w_1 - w_2) \ldots (w_{i-1} - w_i)\\
\qquad \cdot (w_i -w_{i+1}) (w_{i+1} - w_{i+2}) \ldots (w_{n-2} - w_{n-1})
\end{array}
\right]
}  \\
&& \qquad \cdot \left[\sum_{i=k}^{n-2} \left(\frac{1}{v-w_{i+1}} - \frac{1}{v-w_i}\right) + \frac{1}{w_{n-1}-v}\right]\\
&=& \frac{1}{(w_1 - w_2) \ldots (w_{n-2} - w_{n-1})} \frac{1}{w_k - v}\\
&=& \psi_w \cdot \frac{1}{x_{\alpha(e)}-x_{\omega(e)}}
\end{eqnarray*}
The computation is similar if $v$ is the origin of $e$.
\end{proof}

\subsection{Value on forests}
One can now compute the value of $N$ on forests. This result is essential in the following sections because we will often make proofs by induction on the cyclomatic number.

\begin{Prop}\label{prop_arbre}
If $T$ is a tree and $F$ a disconnected forest, one has:
\begin{eqnarray}
N(T)&=&1 ;\\
N(F)&=&0.
\end{eqnarray}
\end{Prop}

\begin{proof}
Thanks to the pruning Lemma \ref{lem:emond} page \pageref{lem:emond}, we only have to prove it in the case where $F$ is a disjoint union of $n$ points. If $n=1$, it is obvious that $N(\cdot) = \Psi(\cdot)=1$. Else, if we denote by $c$ the full cycle $(1 \ldots n)$, one has:
\begin{eqnarray*}
\Psi(F)&=&\sum_{\sigma \in S(n)} \frac{1}{(x_{\sigma(1)}-x_{\sigma(2)}) \ldots (x_{\sigma(n-1)}-x_{\sigma(n)})}\\
&=&\frac{1}{n} \sum_{\sigma \in S(n)} \sum_{i=0}^{n-1} \frac{1}{(x_{\sigma \circ c^i(1)}-x_{\sigma \circ c^i(2)}) \ldots (x_{\sigma \circ c^i(n-1)}-x_{\sigma \circ c^i(n)})}\\
&=&\frac{1}{n} \sum_{\sigma \in S(n)} \frac{\sum\limits_{i=0}^{n-1} x_{\sigma \circ c^i(n)} - x_{\sigma \circ c^i (1)}}{(x_{\sigma(1)}-x_{\sigma(2)}) \ldots (x_{\sigma(n-1)}-x_{\sigma(n)})(x_{\sigma(n)} - x_{\sigma(1)})}\\
&=&0.
\end{eqnarray*}
\end{proof}

\section{The main transformation}\label{sect_oper}
In the section \ref{sect:graphs}, we have defined a simple operation on graphs consisting in removing edges. Thanks to this operation, we will be able to construct an operator which lets invariant the formal sum of linear extensions (paragraph \ref{equ_linea_ext}). Due to the definition of $\Psi$,  this implies immediately an inductive relation on the rational functions $\Psi_G$ (paragraph \ref{cons_rat_func}).\\

\subsection{Equality on linear extensions}\label{equ_linea_ext}
In this paragraph, we prove an induction relation on the formal sums of linear extensions of graphs. More exactly, we write, for any graph $G$ with at least one cycle, $\sumext(G)$ as a linear combination of $\sumext(G')$, where $G'$ runs over graphs with a strictly lower cyclomatic number. In the next paragraphs, we will iterate this relation and apply $\Psi$ to both sides of the equality to study $\Psi_G$.\\

If $G$ is a finite graph and $C$ a cycle of $G$, let us denote by $T_C(G)$ the following formal alternate sum of subgraphs of $G$:
$$T_C(G)=\sum_{\substack{E' \subset HE(C)\\E' \neq \emptyset}} (-1)^{|E'|-1} G \backslash E'.$$

The function $\sumext(G)=\sum\limits_{w \in \ext(G)} w$ can be extended by linearity to the free abelian group spanned by graphs. One has the following theorem:
\begin{Th}\label{th:boucle}
Let $G$ be a graph and $C$ a cycle of $G$. Then,
\begin{equation}\label{eq_boucle}
\sumext(G)=\sumext(T_C(G)).
\end{equation}
\end{Th}

\begin{figure}
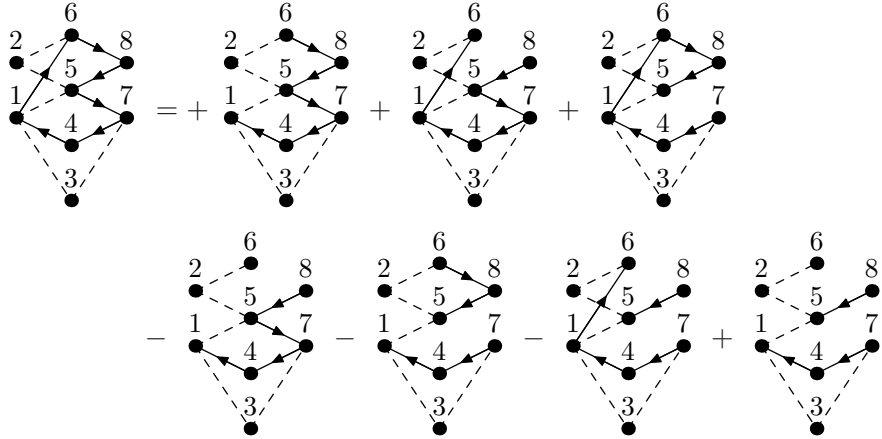


\begin{multline*}
\begin{array}{c}
\mbox{\epsfig{file=figure-8.ps,width=48 pt}}
\end{array}
=
+
\begin{array}{c}
\mbox{\epsfig{file=figure-15.ps,width=48 pt}}
\end{array}
+
\begin{array}{c}
\mbox{\epsfig{file=figure-16.ps,width=48 pt}}
\end{array}
+
\begin{array}{c}
\mbox{\epsfig{file=figure-17.ps,width=48 pt}}
\end{array}\\
-
\begin{array}{c}
\mbox{\epsfig{file=figure-18.ps,width=48 pt}}
\end{array}
-
\begin{array}{c}
\mbox{\epsfig{file=figure-19.ps,width=48 pt}}
\end{array}
-
\begin{array}{c}
\mbox{\epsfig{file=figure-20.ps,width=48 pt}}
\end{array}
+
\begin{array}{c}
\mbox{\epsfig{file=figure-21.ps,width=48 pt}}
\end{array}
\end{multline*}

\caption{Example of application of theorem \ref{th:boucle}}
\label{fig_ex_op}
\end{figure}

Note that all graphs appearing in the right-hand side of (\ref{eq_boucle}) have strictly less cycles than $G$.
An example is drawn on figure \ref{fig_ex_op} (to make it easier to read, we did not write the operator $\sumext$ in front of each graph).\\

\begin{Rem}
In the case where $HE({C})=\emptyset$, this theorem says that graphs with oriented circuits have no linear extensions (see remark \ref{rem_circ} page \pageref{rem_circ}).\\

If it is a singleton, it says that we do not change the set of linear extensions by erasing an edge if there is a path going from its origin to its end (thanks to transitivity).\\

An other very interesting case of our relation is the following one. Let $G$ be a graph and $v_1$ and $v_2$ two vertices of $G$ which are not linked by an edge. We can write
\begin{equation}\label{eq:ord}
\sumext(G)=\sum_{\substack{w \in \ext(G) \\ v_1 \leq_w v_2}} w + \sum_{\substack{w \in \ext(G) \\ v_2 \leq_w v_1}} w 
\end{equation}
This is in fact a special case of our relation on the graph $G'$ obtained from $G$ by adding two edges $e_{1,2}=(v_1,v_2)$ and $e_{2,1}=(v_2,v_1)$. This graph contains a circuit so $\sumext(G')=0$. But one also has:
$$\sumext (G')=\sumext\big(G \cup \{e_{1,2}\} \big) + \sumext\big(G \cup \{e_{2,1}\} \big) - \sumext\big(G \big).$$
So $ \sumext\big(G \big)$ is the sum of two terms corresponding exactly to equation (\ref{eq:ord}). By iterating this equality, deleting graph with circuits and erasing edges thanks to transitivity relation, we obtain:
$$\sumext(G) = \sum_{w \in \ext(G)} \sumext(\begin{array}{c}
                                             \epsfig{file=figure-22.ps,width=3cm}
                                            \end{array}
).$$
An immediate consequence is that any relation between the $\sumext(G)$ can be deduced from Theorem \ref{th:boucle}.
\end{Rem}

To prove Theorem \ref{th:boucle}, we will need the two following lemma:

\begin{Lem}\label{lem:preuve_boucle_1}
Let $w \in \ext(G \backslash HE(C))$. There exists $E'(w)$ such that
$$\forall E'' \subset HE(C), \ \ \ \ w \in \ext(G \backslash E'') \Longleftrightarrow E'(w) \subset E'' \subset HE(C).$$
\end{Lem}

\begin{proof}
immediate consequence of lemma \ref{lem_ext} page \pageref{lem_ext}.
\end{proof}

\begin{Lem}\label{lem:preuve_boucle_2}
Let $w \in \ext(G \backslash HE({C}))$, there exists $E'' \subsetneq HE({C})$ such that
$$w \in \ext(G \backslash E'').$$
\end{Lem}

\begin{proof}
Suppose that we can find a word $w$ for which the lemma is false. Since $w \in \ext(G \backslash HE({C}))$, the word $w$ fulfills the relations of the edges of ${C}$, which are not in $HE({C})$.\\
But, if $e \in HE({C})$, one has $w \notin \ext(G \backslash (HE({C}) \backslash \{e\}))$. That means that $w$ does not fulfill the relation corresponding to the edge $e$. As $w$ is a total order, it fulfills the opposite relation:
$$w \in \ext \left[ \big(G \backslash HE({C})\big) \cup \overline{e}\right].$$
Doing the same argument for each $e \in HE({C})$, one has
$$w \in \ext \left[ \big(G \backslash HE({C})\big) \cup \overline{HE({C})}\right].$$
But this graph contains an oriented cycle so the corresponding set of linear extension is empty.
\end{proof}

Let us come back to the proof of Theorem \ref{th:boucle}.\\
Let $w$ be a word containing exactly once each element of $V(G)$. We will compute its coefficient in $\sumext(G)-\sumext(T_C(G)) = \sum_{E' \subset HE(C)} (-1)^{|E'|} \sumext(G \backslash E')$:
\begin{itemize}
\item If $w \notin \ext(G \backslash HE({C}))$, its coefficient is zero in each summand.
\item If $w \in \ext(G \backslash HE({C}))$, thanks Lemma \ref{lem:preuve_boucle_1}, we know that there exists $E'(w) \subset HE({C})$ such that
$$\forall E'' \subset HE(C), \ \ \ \ w \in \ext(G \backslash E'') \Longleftrightarrow E'(w) \subset E'' \subset HE({C}).$$
So the coefficient of $w$ in $\sumext(G)-\sumext(T_C(G))$ is
$$\sum_{E'(w) \subset E'' \subset HE({C})} (-1)^{|E''|}=0 \text{ (because } E'(w) \neq HE(C) \text{, Lemma \ref{lem:preuve_boucle_2})}.$$
\end{itemize}

\subsection{Consequences on Greene's functions}\label{cons_rat_func}

In the previous paragraph, we have established an induction formula for the formal sum of linear extensions (Theorem \ref{th:boucle}). One can apply $\Psi$ to both sides of this equality to compute $N(G)$:

\begin{Prop}\label{th:boucle_N}
Let $G$ be the graph containing a cycle ${C}$. Then,
$$N(G)=\sum_{\substack{E' \subset HE({C})\\E' \neq \emptyset}} \left[(-1)^{|E'|-1} N(G \backslash E') \prod_{e \in E'} ( x_{\alpha(e)} - x_{\omega(e)} ) \right].$$
\end{Prop}

By Proposition \ref{prop_arbre} page \pageref{prop_arbre}, one has $N(T)=1$ if $T$ is a tree and $N(F)=0$ if $F$ is a disconnected forest. So this Proposition gives us an algorithm to compute $N(G)$: we just have to iterate it with any cycles until all the graphs in the right hand side are forests. More precisely, if after iterating transformations of type $T_C$ on $G$, we obtain the formal linear combination $\sum c_F F$ of subforests of $G$, then:
$$N(G)=\sum_{T \text{ subtree of }G} c_T \prod_{e \in E_G \setminus E_T} (x_{\alpha(e)}-x_{\omega(e)}).$$
In this formula, $N(G)$ appears as a sum of polynomials. So the computation of $N(G)$, using this formula, is easier than a direct application of the definition
$$N(G) = \sum_{w \in \ext(G)} \left( \Psi_w\cdot \prod_{e \in E_G} (x_{\alpha(e)}-x_{\omega(e)}) \right) ,$$
where the summands may have poles.\\

For instance,
$$
N\left(
\begin{array}{c} \mbox{\epsfig{file=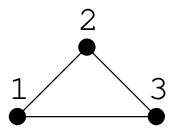,width=48 pt}} \end{array}
\right)
=
N\left(
\begin{array}{c} \mbox{\epsfig{file=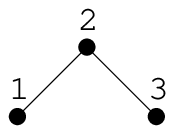,width=48 pt}} \end{array}
\right) . (x_1-x_3) = x_1 -x_3.
$$
We will use this algorithm in the next section on some other examples. But it has also a theoretical interest: some properties of $N$ on forests can be immediately extended to any graph.

\begin{Corol}\label{Corol:zero_nc}
For any graph $G$, the rational function $N(G)$ is a polynomial. Moreover, if $G$ is disconnected, $N(G)=0$.
\end{Corol}

In fact, if $G$ is the Hasse diagram of a connected poset, the fraction $\Psi(G)=\frac{N(G)}{\prod\limits_{e\in E_G} (x_{\alpha(e)}-x_{\omega(e)})}$ is irreducible (see \cite{Bo} for a proof of this fact).

The following result can also be proved by induction on the cyclomatic number:

\begin{Prop}
\label{prop_contraction}
Let $G$ be a graph and $e$ an edge of $G$ between two vertices $v_1$ and $v_2$. Then
$$N(G/e) = N(G) \big|_{x_{v_1}=x_{v_2}=x_v},$$
where $v$ is the contraction of $v_1$ and $v_2$ in $G/e$.
\end{Prop}

\begin{figure}
$$
\begin{CD}
\begin{array}{c} \mbox{\epsfig{file=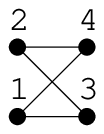,width=24 pt}} \end{array}  @>\text{contraction of }2-4>> \begin{array}{c} \mbox{\epsfig{file=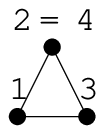,width=24 pt}} \end{array}\\
@VNVV  @VVNV\\
x_1+x_2-x_3-x_4  @>x2=x4>> x1-x3
\end{CD}
$$
\caption{Contraction and numerator}
\label{fig_example_contraction}
\end{figure}

\begin{proof}[Proof (by induction on the cyclomatic number of $G$)]
If $G$ is a forest, then the equality is obvious by Proposition \ref{prop_arbre}.\\

If $G/e$ contains a cycle $C_e$, then we consider the following cycle $C$ in $G$ (figure \ref{fig_case_proof_contraction} illustrate all the different cases):

\begin{figure}
$$
\begin{CD}
G = \! \! \! \! \begin{array}{c} \mbox{\epsfig{file=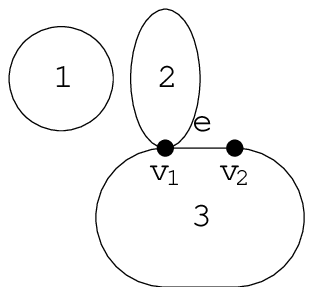,width=100 pt}} \end{array}  @>\text{contraction} >\text{ of }e=(v_1,v_2)> G / e = \! \! \! \! \begin{array}{c} \mbox{\epsfig{file=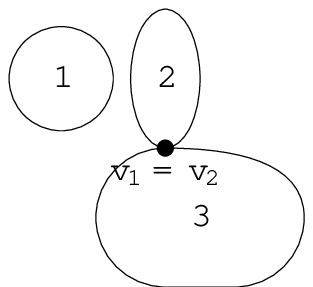,width=100 pt}} \end{array}
\end{CD}
$$
\caption{The different cases of the proof of proposition \ref{prop_contraction}}
\label{fig_case_proof_contraction}
\end{figure}

\begin{enumerate}[1)]
 \item If $C_e$ does not go through the vertex $v$ (contraction of $v_1$ and $v_2$), then $C_e$ can also be seen as a cycle $C$ of $G$.
 \item Suppose that $v$ is the end of $e_i$ and the origin of $e_{i+1}$ and that they are also the same vertex ($v_1$ or $v_2$) in $G$. Then, $C_e$ can still be seen as a cycle $C$ of $G$.
 \item Suppose that $v$ is the end of $e_i$ and the origin of $e_{i+1}$ but that these two edges have different extremities ($v_1$ and $v_2$) in $G$. Then we add the edge $e$ or $\overline{e}$ to $C_e$ (between $e_i$ and $e_{i+1}$) to obtain a cycle $C$ of $G$.
\end{enumerate}

Eventually by changing the orientations of $C_e$ and $C$, we can assume that $e \notin HE(C)$ and, as a consequence $HE(C)=HE(C_e)$. By theorem \ref{th:boucle_N} page \pageref{th:boucle_N}, one has:
\begin{eqnarray*}
 N(G/e)=\sum_{\substack{E' \subset HE(C_e)\\E' \neq \emptyset}} (-1)^{|E'|-1} N( (G/e) \backslash E').\prod_{e \in E'}(x_{\alpha(e)} - x_{\omega(e)} ) \\
 N(G)=\sum_{\substack{E' \subset HE(C)\\E' \neq \emptyset}} (-1)^{|E'|-1}N(G \backslash E').\prod_{e \in E'}(x_{\alpha(e)} - x_{\omega(e)}).
\end{eqnarray*}
As $e \notin HE(C)$,
$$(G \backslash E')/e=(G/e) \backslash E' \ \text{ and } \  HE(C_e) = HE(C).$$
This ends the proof by applying the induction hypothesis to the graphs $G \backslash E'$.
\end{proof}

Another immediate consequence of Proposition \ref{th:boucle_N} is the following vanishing property of $N(G)$.
\begin{Corol}\label{corol:annulation}
Let $G$ be a graph. Let $C$ be a cycle of $G$ with $HE({C})=\{e_1,\ldots,e_r\}$. One has
$$\left.N(G)\right|_{ x_{\alpha(e_i)}=x_{\omega(e_i)}, i=1 \ldots r}=0$$
\end{Corol}
Unfortunately, this corollary, written for every cycle of a graph $G$, does not characterize $N(G)$ up to a multiplicative factor (see paragraph \ref{subsubsect:openpb_annulation}).

\section{Some explicit computations of rational functions} \label{sect_exemples}
This section is devoted to some examples of explicit computation of $N(G)$ using the algorithm described in paragraph \ref{cons_rat_func}.

\subsection{Graphs with cyclomatic number 1.}\label{subsect:ex_1cycle}
We consider in this paragraph connected graphs $G$ with $|V_G|=|E_G|$. Using pruning Lemma \ref{lem:emond} page \pageref{lem:emond}, we can suppose that each vertex of $G$ has valence $2$. We denote by $\max(G)$ (resp. $\min(G)$) the set of maximal (resp. minimal) elements of $G$. The following result was already proved in \cite{Bo}, but we present here a simpler proof using the results of the previous section.

\begin{Prop}
If $G$ is a connected graph with vertices of valence $2$, then
$$N(G) = \sum_{v \in \min(G)} x_v - \sum_{v' \in \max(G)} x_{v'}.$$
\end{Prop}

\begin{proof}
 $G$ has only one cycle $C$ (we only have to choose an orientation). In the right-hand side of equation (\ref{eq_boucle}) page \pageref{eq_boucle}, we have two kinds of terms:
\begin{itemize}
 \item If $|E'|=1$, $G \backslash E'$ is a tree and $N(G \backslash E')=1$.
 \item If $|E'|>1$, $G \backslash E'$ is disconnected and $N(G \backslash E')=0$.
\end{itemize}
Then
$$N(G)=\sum_{e \in HE(C)} ( x_{\alpha(e)} - x_{\omega(e)}).$$
The sum above can be simplified and is equal to $\sum\limits_{v \in \min(G)} x_v - \sum\limits_{v' \in \max(G)} x_{v'}$.
\end{proof}

\begin{example}
\begin{eqnarray*}
N \left(\begin{array}{c}
\mbox{\epsfig{file=figure-29.ps,width=30 pt}}
\end{array}\right)
&=&
(x_1-x_3) N \left(\begin{array}{c}
\mbox{\epsfig{file=figure-30.ps,width=30 pt}}
\end{array}\right)
+
(x_2 - x_4) N \left(\begin{array}{c}
\mbox{\epsfig{file=figure-31.ps,width=30 pt}}
\end{array}\right)\\
& & +
(x_4-x_5)N \left(\begin{array}{c}
\mbox{\epsfig{file=figure-32.ps,width=30 pt}}
\end{array}\right)
\pm
N \left(\begin{array}{c} \text{disconnected}\\ \text{graphs} \end{array} \right)\\
&=& (x_1-x_3) + (x_2 - x_4) + (x_4-x_5)\\
&=& x_1 + x_2 - x_3 - x_5.
\end{eqnarray*}
\end{example}

\subsection{Graphs with cyclomatic number 2.} \label{subsect:graphe2boucles}
Let $G$ be a connected graph with a cyclomatic number equal to $2$. Thanks to pruning Lemma \ref{lem:emond} page \pageref{lem:emond}, we can assume that $G$ has no vertices of valence $1$. As $|E_G|=|V_G|+1$, the graph has, in addition of vertices of valence $2$, either two vertices of valence $3$ or one vertex of valence $4$. We will only look here at the case where there are two vertices $v$ and $v'$ of valence $3$ and the edges can be partitioned into three paths $p_0$, $p_1$ and $p_2$ from $v$ to $v'$ (the other cases are easier because the cycles have no edges in common). \\

\begin{figure}
\begin{center}
\mbox{\epsfig{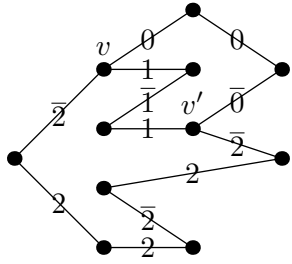}}
\caption{Example of a graph $G$ with cyclomatic number $2$.}
\label{fig_2_boucles}
\end{center}
\end{figure}

For $i=0,1,2$, let us denote by $E_i$ (resp. by $E_{\overline{i}}$) the set of edges of the path $p_i$ which appear in the same (resp. opposite) orientation in the graph and in the path $p_i$ (see the figure \ref{fig_2_boucles} for an example, we have written on each edge the index of the set it belongs to). If $I=\{i_1,\ldots,i_l\} \subset \{0,1,2,\overline{0},{\overline{1}},{\overline{2}}\}$, we consider the following alternate sum of graphs:
\begin{equation}
\label{eq:graphs_2}
G_I = \sum_{\emptyset \neq E'_1 \subset E_{i_1} , \ldots , \emptyset \neq E'_l \subset E_{i_l}} (-1)^{|E'_1|-1} \ldots (-1)^{|E'_l|-1} G \backslash \big( E'_1 \cup \ldots \cup E'_l \big).
\end{equation}
Let us consider the cycle ${C}=\overline{p_1} \cdot p_2$: one has $HE(C)= E_{\overline{1}} \cup E_2$. The subsets of $HE(C)$ can be partitioned in three families:
\begin{itemize}
\item that included in $E_{\overline{1}}$ ;
\item that included in $E_2$;
\item the unions of non empty subset of $E_{\overline{1}}$ and non empty subset of $E_2$.
\end{itemize}
Thus, if we apply Theorem \ref{th:boucle} page \pageref{th:boucle} with respect to $C$, we obtain:
$$\sumext(G) = \sumext(G_{\overline{1}}) + \sumext(G_2) - \sumext(G_{2,{\overline{1}}}),$$
where $G_1$,$G_2$ and $G_{2,1}$ are defined in equation \ref{eq:graphs_2}.

Each graph in $G_{\overline{1}}$ contains the cycle $\overline{p_0} \cdot p_2$, because only edges belonging to $p_1$ have been removed. If we apply Theorem \ref{th:boucle} with this cycle, we obtain:
\begin{multline*}
\sumext(G_{\overline{1}}) = \sum_{E' \subset E_{\overline{1}}} (-1)^{|E'|-1}\sumext(G \setminus E') \\
= \sum_{E' \subset E_{\overline{1}}} (-1)^{|E'|-1} \left( \sum_{E''\subset E_{\overline{0}}} (-1)^{|E''|-1} \sumext\big((G \setminus E') \setminus E'' \big) \right. \\ +  \sum_{E''\subset E_{2}} (-1)^{|E''|-1} \sumext\big((G \setminus E') \setminus E'' \big) \\ \left. -  \sum_{\substack{E''\subset E_{2} \\ E'''\subset E_{\overline{0}}}} (-1)^{|E''|-1} (-1)^{|E'''|-1} \sumext\big((G \setminus E') \setminus (E'' \cup E''') \big) \right)\\
= \sumext(G_{{\overline{0}},{\overline{1}}}) + \sumext(G_{2,{\overline{1}}}) - \sumext(G_{2,{\overline{0}},{\overline{1}}})
\end{multline*}
In a similar way, all graphs in $G_2$ contains the cycle $p_1 \cdot \overline{p_0}$ and one has $\sumext(G_2) = + \sumext(G_{2,{\overline{0}}}) + \sumext(G_{1,2}) -  \sumext(G_{1,2,{\overline{0}}})$. The graphs in $G_{2,{\overline{1}}}$ have no cycles, so, finally:
\begin{eqnarray*}
\sumext(G)&=&  \sumext(G_{{\overline{0}},{\overline{1}}}) + \sumext(G_{2,{\overline{1}}}) - \sumext(G_{2,{\overline{0}},{\overline{1}}}) \\
& & \qquad + \sumext(G_{2,{\overline{0}}}) + \sumext(G_{1,2}) -  \sumext(G_{1,2,{\overline{0}}}) -\sumext(G_{2,{\overline{1}}}) ;\\
&=&  \sumext(G_{{\overline{0}},{\overline{1}}}) - \sumext(G_{2,{\overline{0}},{\overline{1}}}) + \sumext(G_{2,{\overline{0}}}) + \sumext(G_{1,2}) - \sumext(G_{1,2,{\overline{0}}}).
\end{eqnarray*}
If we apply $\Psi$ to this equality, we keep only connected graphs and obtain:
$$\Psi(G) = \Psi(G'_{{\overline{0}},{\overline{1}}}) + \Psi(G'_{2,{\overline{0}}}) + \Psi(G'_{1,2}),$$
where $G'_I= \sum_{e_1 \in E_{i_1} , \ldots , e_l \in E_{i_l}} G \backslash \big\{ e_1 , \ldots , e_l \big\}.$
As all graphs in the expression of $G'_I$ are trees, we obtain ( by using $X_e$ instead of $x_{\alpha(e)} - x_{\omega(e)}$ ):
\begin{eqnarray*}
N(G) & = &  \sum_{e_{\overline{0}}  \in E_{\overline{0}}, e_{\overline{1}} \in E_{\overline{1}}} \! X_{e_{\overline{0}}} \cdot X_{e_{\overline{1}}}
+ \! \sum_{e_{\overline{0}} \in E_{\overline{0}}, e_2 \in E_2} \! X_{e_{\overline{0}}} \cdot X_{e_2}
+ \! \sum_{e_1 \in E_1, e_2 \in E_2} \!  X_{e_1} \cdot X_{e_2}\\
& = & \left( \sum_{e_{\overline{0}} \in E_{\overline{0}}} X_{e_{\overline{0}}} \right) \left( \sum_{e_{\overline{1}} \in E_{\overline{1}}} X_{e_{\overline{1}}} \right) +
\left( \sum_{e_{\overline{0}} \in E_{\overline{0}}} X_{e_{\overline{0}}} \right) \left( \sum_{e_2 \in E_2} X_{e_2} \right) \\
& & \qquad + \left( \sum_{e_1 \in E_1} X_{e_1} \right) \left( \sum_{e_2 \in E_2} X_{e_2} \right).
\end{eqnarray*}
One can notice that, if $E_{\overline{0}}$ is empty (that is to say that there is a chain form $v$ to $v'$), the polynomial $N(G)$ is the product of two polynomials on degree $1$. This is a particular case of our third main result (Theorem \ref{factorization_chaine}).

\subsection{Simple bipartite graphs} \label{subsect:bipartite}
\begin{Def}\label{def_bipartite}
A graph is said to be bipartite if its set of vertices can be partitioned in two sets $A= \{ a_i \}$ and $B = \{ b_i \}$ such that $E \subset A \times B$.\\
Moreover, a bipartite graph is said complete if $E=A \times B$.
\end{Def}

\begin{figure}
\begin{center}
\mbox{\epsfig{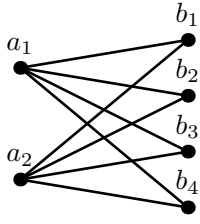}}
\caption{The bipartite graph $G_{2,4}$.}
\label{fig_bipartite}
\end{center}
\end{figure}

In this section we will look at bipartite graphs $G$ such that $|A|=2$. Thanks to the pruning Lemma \ref{lem:emond} page \pageref{lem:emond}, we can suppose that $G$ is a complete bipartite graph. The complete bipartite graph with $|A|=2$ and $|B|=n$ is unique up to isomorphism and will be denoted $G_{2,n}$ (drawn on figure \ref{fig_bipartite} for $n=4$).\\

We will denote vertices and his associated variables in the same way.

\begin{Prop}
Let $G_{2,n}$ be a bipartite graph with $A=\{a_1,a_2\}$ and $B = \{ b_1,\ldots,b_n \}$, then
$$N(G_{2,n})= \sum_{i=1}^n \left( \prod_{j < i} (a_1 - b_j) \cdot \prod_{k > i} (a_2 - b_k) \right).$$
\end{Prop}

\begin{proof}
For each $h=1,2$ and $i=1,\ldots,n$, we denote by $e_{h,i}$ the edge $(a_h,b_i)$.
We will show, by induction of $n$, that, by applying several times theorem \ref{th:boucle} page \pageref{th:boucle}, we obtain the following equality (which is drawn on figure \ref{fig_G24} for $n=4$ ; we omit the $\sumext$ for clearness):
\begin{multline}
\sumext(G_{2,n}) = \sum_{i=1}^n \sumext\big(G_{2,n} \backslash \big\{e_{2,1},\ldots,e_{2,i-1},e_{1,i+1},\ldots,e_{1,n}\big\}\big)\\
 - \sum_{i=1}^{n-1} \sumext \big( G_{2,n} \backslash \{e_{2,1},\ldots,e_{2,i},e_{1,i+1},\ldots,e_{1,n} \big\}\big).
\end{multline}

\begin{figure}
\begin{center}
$
\begin{array}{c}
\mbox{\epsfig{file=figure-35.ps,width=24 pt}}
\end{array}
=
\begin{array}{c}  
\mbox{\epsfig{file=figure-36.ps,width=24 pt}}
\end{array}
+
\begin{array}{c}
\mbox{\epsfig{file=figure-37.ps,width=24 pt}}
\end{array}
+
\begin{array}{c}
\mbox{\epsfig{file=figure-38.ps,width=24 pt}}
\end{array}
+
\begin{array}{c}
\mbox{\epsfig{file=figure-39.ps,width=24 pt}}
\end{array}
-
\begin{array}{c}

\mbox{\epsfig{file=figure-40.ps,width=24 pt}}
\end{array}
-
\begin{array}{c}
\mbox{\epsfig{file=figure-41.ps,width=24 pt}}
\end{array}
-
\begin{array}{c}
\mbox{\epsfig{file=figure-42.ps,width=24 pt}}
\end{array}
$
\end{center}
\caption{Decomposition of $\sumext(G_{2,4})$.\label{fig_G24}}
\end{figure}

For $n=1$, the statement is obvious. Let us suppose that our formula is true for $n$ and that the equality at rank $n$ is obtained by an iterated application of Theorem $\ref{th:boucle}$ in the graph $G_{2,n}$. We can do the same transformations in $G_{2,n+1}$ (which contains canonically $G_{2,n}$). We obtain:
\begin{multline}\label{eqtechquad}
\sumext(G_{2,n+1}) = \sum_{i=1}^n \sumext\big(G_{2,n+1} \backslash \big\{e_{2,1},\ldots,e_{2,i-1},e_{1,i+1},\ldots,e_{1,n}\big\}\big) \\
- \sum_{i=1}^{n-1} \sumext \big( G_{2,n+1} \backslash \{e_{2,1},\ldots,e_{2,i},e_{1,i+1},\ldots,e_{1,n} \big\}\big)
\end{multline}
The graphs of the first line have still one cycle ($e_{2,i},\overline{e_{1,i}},e_{1,n+1},\overline{e_{2,n+1}}$). By Theorem \ref{th:boucle}, one has: 
\begin{multline*}\sumext(G_{2,n+1} \backslash \{e_{2,1},\ldots,e_{2,i-1},e_{1,i+1},\ldots,e_{1,n}\}) =\\
 \sumext \big( G_{2,n+1} \backslash \{e_{2,1},\ldots,e_{2,i-1},e_{2,i},e_{1,i+1},\ldots,e_{1,n}\}\big)\\
 + \sumext \big( G_{2,n+1} \backslash \{e_{2,1},\ldots,e_{2,i-1},e_{1,i+1},\ldots,e_{1,n},e_{1,n+1}\} \big)\\
 - \sumext \big( G_{2,n+1} \backslash \{e_{2,1},\ldots,e_{2,i-1},e_{2,i},e_{1,i+1},\ldots,e_{1,n},e_{1,n+1}\big).
\end{multline*}
Using this formula for each $i$, the first summand balances with the negative term in (\ref{eqtechquad}) (except for $i=n$) and the two other summands are exactly what we wanted. This ends the induction and Formula (\ref{eqtechquad}) is true for any $n$.\\

Note that the graphs of its right hand side have no cycles and that only the ones of the first line are connected. We just have to apply $\Psi$ to this equality, and use the value of $\Psi$ on forests (Proposition \ref{prop_arbre} page \pageref{prop_arbre}) to finish the proof of the proposition.
\end{proof}

Note that this case is interesting because the function $N$ can be expressed as a specialization of a rectangular Schur function (see \cite[Proposition 4.2]{Bo}).

\begin{Rem}
Our algorithm allows us to write $\sumext(G)$ as a sum of terms of the kind $\pm \sumext(F)$, with $F$ subforest of $G$. But, in the three examples of this section, all trees have $0$ or $+1$ as coefficients after iteration of transformations of kind $T_C$ on $G$. We will see in the next section that this is possible for any graph $G$ with a clever choice of cycles.
\end{Rem}

\section{A combinatorial formula for N} \label{sect:combi}
To compute the polynomial $N$ associated to a graph $G$, we only have to find the coefficient of trees in a formal linear combination of forests obtained by iterating transformations $T_C$ on $G$. But there are many possible choices of cycles at each step and these coefficients depend on these choices.\\

A way to avoid this problem is to give to $G$ a rooted map structure and to look at the particular decomposition introduced in the paper \cite[section 3]{Fe}. With these particular choices, we have a combinatorial description of the trees with coefficient $+1$, all other trees having $0$ as coefficient.

\subsection{Rooted maps and admissible cycles}

\begin{Def}
A (combinatorial oriented) map is a connected graph with, for each vertex $v$, a cyclic order on the edges whose origin or end is $v$. This definition is natural when the graph is drawn on a two dimensional surface (see for example \cite{CARTES}). The figure \ref{fig_example_map} gives an example of three different maps with the same underlying graph. \\

\begin{figure}
\begin{center}
\epsfig{file=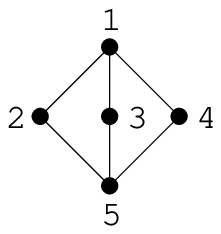,width=60 pt}
\hspace{1 cm}
\epsfig{file=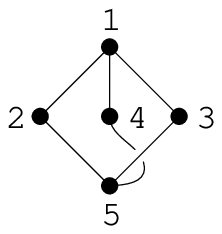,width=60 pt}
\hspace{1 cm}
\epsfig{file=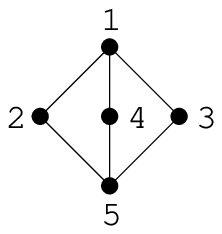,width=60 pt}
\end{center}
\caption{Example of three different maps.\label{fig_example_map}}
\end{figure}

It is more convenient when we deal with maps, to consider edges as couples of two half-edges (called darts) $(h_1,h_2)$, the first one of extremity $\alpha(e)$ and the second one of extremity $\omega(e)$. Then the map structure is given by a permutation $\sigma$ of all the darts, whose orbits correspond to the sets of darts with the same extremity.\\

A rooted map is a map with an external dart $h_0$, that is to say a dart which do not belong to any edge, but has an extremity (which will be denoted by $\star$) and a place in the cyclic order given by this extremity.\end{Def}

\begin{Rem}
In this section, as cyclic orders of edges around vertices matter, we can not use the convention that the extremity of an edge is always on its origin's right (we did not assume any condition on compatibility between the orientations of the edges and the map structure, see open problem \ref{subsubsect:openpb_mapstructure}).
\end{Rem}

Recall that, to compute $N(G)$, a naive algorithm is to choose any cycle of the graph, apply proposition \ref{th:boucle_N} page \pageref{th:boucle_N}. If the graph has a rooted map structure, it is interesting to choose cycles with additional properties. Our choices will not involve the orientation of the edges of the map. So we will define a notion of admissible cycle in a (not necessary oriented) rooted map.\\

\begin{figure}
\begin{center}
$$\epsfig{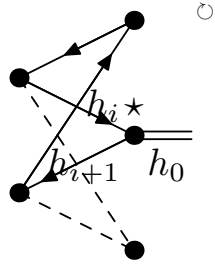}^\circlearrowright$$
\caption{Example of a rooted map $M$ with an admissible cycle of type 1.}
\label{fig:ex_adm1}
\end{center}
\end{figure}

By definition, a cycle ${C}$ of a rooted map is admissible of type $1$ (see figure \ref{fig:ex_adm1}) if:
\begin{itemize}
\item The vertex $\star$ is a vertex of the cycle, that is to say that $\star$ is the extremity of a dart $h_{i}$ of $e_i$ and of a dart $h_{i+1}$ of $e_{i+1}$ for some $i$ ;
\item The cyclic order at $\star$ restricted to the set $\big\{h_0,h_{i},h_{i+1}\big\}$ is the cyclic order $\big(h_0,h_{i+1},h_{i} \big)$.
\end{itemize}

If ${C}$ satisfies the first condition, exactly one cycle among ${C}$ and ${\overline{C}}$ is admissible (where ${\overline{C}}$ is ${C}$ with the opposite orientation).\\

\begin{figure}
\begin{center}
\includegraphics{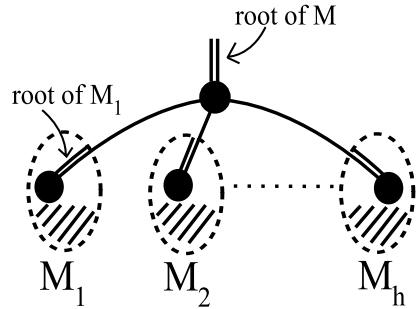}\caption{A generic rooted map $M$ without admissible cycles of type $1$}
\label{figsortieD1}
\end{center}
\end{figure}

\begin{figure}
\begin{center}
$$\epsfig{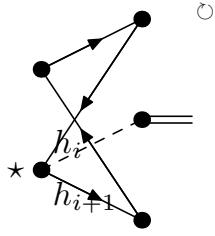}^\circlearrowright$$
\caption{Example of a rooted map $M'$ with an admissible cycle of type 2.}
\label{fig:ex_adm2}
\end{center}
\end{figure}

If a rooted map has no admissible cycles of type $1$, it is of the form of the figure \ref{figsortieD1}. In this case, we call admissible cycles of type $2$ the admissible cycles of its "legs" $M_1,\ldots,M_h$ (of type $1$ or $2$, this defines the admissible cycles by induction). Note that this definition has a sense because the legs have a canonical external dart and are rooted maps. An example of an admissible cycle of type $2$ is drawn on Figure \ref{fig:ex_adm2}\\

A rooted map without admissible cycles has no cycles at all, hence it is a tree.

\begin{Rem}
 The second condition in the definition of admissible of type $1$ says that the root must be at the left of the cycle. The first condition is only technical, because if the cycle does not go through $\star$, we can not define \emph{``to be on the left of the cycle''}.\\
For a planar map this can be avoided because any cycle split the plan into two regions, so the left side of an oriented cycle is well-defined. In this case, we can call admissible any cycle such that the root is at the left of the cycle even if the the cycle does not go through $\star$ and the confluence of the algorithm in the next paragraph will still be true.
\end{Rem}

\subsection{Decomposition of rooted maps}

Consider the following algorithm:
\begin{description}
\item[Input] a rooted map $M$.
\item[Variable] $S$ is a formal linear combination of submap of $M$.
\item[Initialization] $S=M.$
\item[Iterated step] Choose a map $M_0$ with a non-zero coefficient $c_{M_0}$ in $S$ which is not a forest and $C$ an admissible cycle of $M_0$. Apply $T_C$ to $M_0$ in $S$ and keep only the connected graphs in the right-hand side (they have a natural induced rooted map structure). Formally,
$$S:= S - c_{M_0} M_0 + c_{M_0} \delta(T_C(M_0)),$$
where $\delta$ is the linear operator defined by:
$$\delta(M')=\left\{ \begin{array}{ll}
                      M' & \text{if $M'$ connected}\\
		      0 & \text{else}
                     \end{array} \right. .$$
\item[End] We iterate this until $S$ is a linear combination of subtrees of $M$.
\item[Output] $S$.
\end{description}

\begin{Defth}
This algorithm always terminates and is confluent. Let $D(M)$ be its output.
\end{Defth}

\begin{proof}[Idea of the proof]
The termination is obvious: all maps in $T_C(M_0)$ have a lower cyclomatic number than $M_0$.\\

For the confluence, the maps whose graphs are considered in paragraph \ref{subsect:graphe2boucles} play a similar role to critical peaks in rewriting theory. We just have to check our result on these maps. There are infinitely many maps of this kind, but, as in paragraph \ref{subsect:graphe2boucles}, one computation is enough to deal with the general case.
\end{proof}
For a complete proof, see \cite[definition-theorem 3.1.1 and 3.2.1, together with remark 2]{Fe}.

\begin{Prop}
 Let $M$ be a rooted map.
$$\Psi(D(M))=\Psi(M)$$
\end{Prop}
\begin{proof}
 We have to check that $\Psi(S)$ is an invariant of our algorithm. This is trivial because operators $T_C$ and $\delta$ let $\Psi$ invariant (see Theorem \ref{th:boucle} page \pageref{th:boucle} and Corollary \ref{Corol:zero_nc} page \pageref{Corol:zero_nc}). 
\end{proof}

\begin{example}
Let $M$ be the complete bipartite graph $G_{2,3}$ ($A=\{a_1,a_2\}$, $B=\{b_1,b_2,b_3\}$) with the following rooted map structure:
\begin{itemize}
 \item If we denote by $e_{1,i}$ (resp $e_{2,i}$) the edge between $a_1$ (resp. $a_2$) and $b_i$, the cyclic order around the vertex $a_1$ (resp. $a_2$) is $(e_{1,1},e_{1,2},e_{1,3})$ (resp. $(e_{2,1},e_{2,2},e_{2,3})$).
 \item The root has extremity $b_2$ and is located before $e_{2,2}$.
\end{itemize}
The cycle $C=(\overline{e_{2,2}},e_{2,1},\overline{e_{1,1}},e_{1,2})$ with $HE(C)=\{e_{2,1},e_{1,2}\}$ (drawn on Figure \ref{fig:ex_adm1}) is admissible (of type 1). So, with this choice, after the first iteration of step $1$ of our decomposition algorithm, we have:

$$S=
\begin{array}{c}
\mbox{\epsfig{file=figure-49.ps,width=40pt}}
\end{array}
+\begin{array}{c}
\mbox{\epsfig{file=figure-50.ps,width=40pt}}
\end{array}
-\begin{array}{c}
\mbox{\epsfig{file=figure-51.ps,width=40pt}}
\end{array}$$
The two firsts graph have each an admissible cycle: the first one of type 1 ($C=(\overline{e_{2,2}},e_{2,3},\overline{e_{1,3}},e_{1,2})$ with $HE(C)=\{e_{2,3},e_{1,2}\}$), the second one of type 2 ($C=(e_{2,3},\overline{e_{1,3}},e_{1,1},\overline{e_{2,1}})$ and $HE(C)=\{e_{2,3},e_{1,1}\}$, see figure \ref{fig:ex_adm2}). So the algorithm ends after two other iterations and we obtain:
\begin{eqnarray}
D(M)&=&\begin{array}{c}
\mbox{\epsfig{file=figure-52.ps,width=30pt}}
\end{array}
+\begin{array}{c}
\mbox{\epsfig{file=figure-51.ps,width=30pt}}
\end{array}
+\begin{array}{c}
\mbox{\epsfig{file=figure-53.ps,width=30pt}}
\end{array}
+\begin{array}{c}
\mbox{\epsfig{file=figure-54.ps,width=30pt}}
\end{array}
-\begin{array}{c}
\mbox{\epsfig{file=figure-51.ps,width=30pt}}
\end{array}\\
&=&\begin{array}{c}
\mbox{\epsfig{file=figure-52.ps,width=30pt}}
\end{array}
+\begin{array}{c}
\mbox{\epsfig{file=figure-53.ps,width=30pt}}
\end{array}
+\begin{array}{c}
\mbox{\epsfig{file=figure-54.ps,width=30pt}}
\end{array}
\end{eqnarray}
Note that, after cancellation, the coefficient of trees in $D(M)$ are $0$ or $+1$. In the next paragraph we will show that it is true for any map $M$ (the sign is a particular case of \cite[Proposition 3.3.1]{Fe}) and characterize combinatorially the trees with a coefficient $+1$.
\end{example}

\subsection{Coefficients in $D(M)$}

To compute the polynomial $N$, we only have to compute the coefficients of spanning trees in $D(M)$. In this section, we will link this coefficient with a combinatorial property of the tree $T$.

\begin{figure}
$$
\epsfig{file=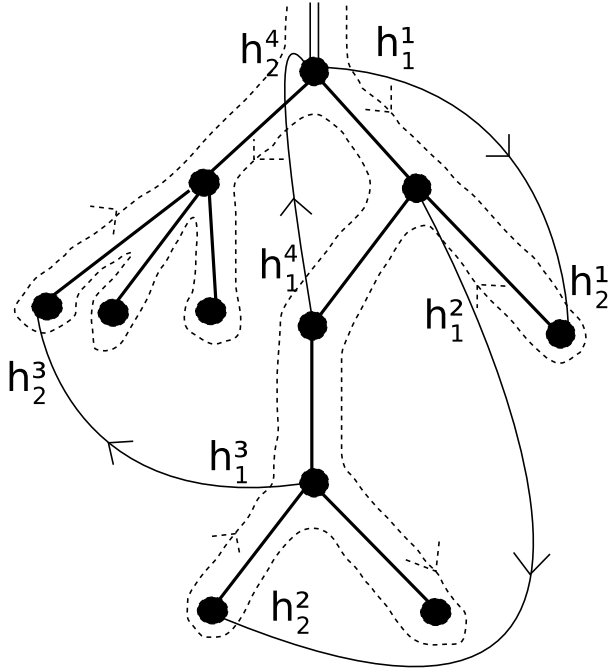,width=8cm}$$
\caption{Tour of a spanning tree of a map.}
\label{fig_tourarbre}
\end{figure}

\begin{Def}
If $T$ is a spanning subtree of a rooted map $M$, the tour of the tree $T$ beginning at $h_0$ defines an order on the darts which do not belong to $T$. The definition is easy to understand on a figure: for example, on Figure \ref{fig_tourarbre}, the tour is ($h^1_1,h^1_2,h^2_1,h^2_2,h^3_1,h^4_1,h^3_2,h^4_2$).  (see \cite{Be} for a precise definition).
\end{Def}

We recall that $D(M)$ does not depend on the admissible cycle chosen at step $1$ of the decomposition algorithm. A good choice to compute the coefficient of a given spanning tree $T \subsetneq M$ is given by lemma \ref{lem_chboucle}. Given an edge $e$ of $M \backslash T$, it is well-known that there exists a unique cycle (up to the orientation) denoted $C(e)$ such that $C(e) \subset (E_T \cup \{e\})$. 

\begin{Lem}\label{lem_chboucle}

There exists an edge $e_0 \in M \backslash T$ such that $C(e_0)$ or $\overline{C(e_0)}$ is admissible. We assume without loss of generality that $C(e_0)$ is admissible. Moreover,
$$e_0 \in HE(C(e_0)) \Longleftrightarrow \begin{array}{c} \text{The first dart of $e_0$ appears in} \\ \text{the tour of $T$ before the second one.}\end{array}$$
\end{Lem}

\begin{proof}
The proof of the lemma, by induction on the size of $M$, can be divided in three cases:
\begin{enumerate}
\item If there is an edge $e$ of $M \backslash T$ whose origin or end is $\star$ (the extremity of the external dart), then $\star$ is a vertex of the cycle $C(e)$ and either $C(e)$ or $\overline{C(e)}$ is admissible of type $1$.
\item Else, let $T_1, \ldots, T_l$ be the connected component of $T \backslash \{\star\}$. If there is an edge $e$ whose extremities are in two different $T_i$, then $C(e)$ is going through $\star$ and $e$ suits in the lemma. 
\item Else, $M \backslash \{\star\}$ has as many connected components as $T \backslash \{\star\}$. Let us denote them by $M_i \supset T_i (1 \leq i \leq l)$. There exists an $j$, such that $M_j \supsetneq T_j$. In this case $M$ has no admissible cycle of type $1$, but by induction there exists $e \in M_j \backslash T_j$ such that $C_{M_j}(e)$ is admissible in $M_j$. By definition, this cycle is admissible of type $2$ in $M$. But $C_{M_j}(e)=C_{M}(e)$, so the proof of the lemma is over.
\end{enumerate}
The second part of the proof is easy  in the two first cases (see figure \ref{fig_tourarbre}). For the third one, it is again an immediate induction.
\end{proof}

This helps us to compute all coefficients of trees in $D(M)$:

\begin{Prop}\label{prop:coef_arbres}
Let $M$ be a rooted map and $T$ a spanning tree of $M$.
\begin{itemize}
\item If there is an edge $e=(h_1,h_2) \in M \backslash T$ such that $h_2$ appears before $h_1$ in the tour of $T$, then the coefficient of $T$ in $D(M)$ is $0$.
\item Else, the coefficient of $T$ in $D(M)$ is $+1$ ($T$ will be said \emph{good}). 
\end{itemize}
\end{Prop}
For example, the spanning tree of Figure \ref{fig_tourarbre} is good. Note that the property of being a good spanning tree does not depend on the orientation of the edges of the tree, but only on the orientation of those which do not belong to it (which is represented by arrows on Figure \ref{fig_tourarbre}.

\begin{proof}
We will prove this proposition by induction over the number of edges in $M \backslash T$. If $M=T$, $T$ is good and the result is obvious.\\

Let $T$ be a covering tree of rooted map $M$ such that $M \backslash T$ contains at least one element. From lemma \ref{lem_chboucle} page \pageref{lem_chboucle}, there exists an edge $e_0$ such that $C(e_0)$ is admissible.
Two cases have to be examinated:
\begin{description}
\item[Case $e_0 \notin HE(C(e_0))$] In this case the spanning tree $T$ can not be good. Besides, $HE(C(e_0)) \subset T$, so every map appearing in $T_{C(e_0)}(M)$ does not contain $T$. But this remains true when we apply operators of kind $T_C$. In particular, the coefficient of $T$ in $D(M)$ is $0$.
\item[Case $e_0 \in HE(C(e_0))$] In this case, one has:
\begin{eqnarray*}T_{C(e_0)}(M) & = & M \backslash \{e_0\} + \text{ maps which do not contain } T.\\
\text{So } D(M) & = & D(M \backslash \{e_0\}) + \sum_{M' \nsupseteq T} D(M').
\end{eqnarray*}
As in the previous case, the second summand has a contribution $0$ to the coefficient of $T$ in $D(M)$. By induction hypothesis, the first one has contribution $+1$ if $T$ is a good spanning tree of $M \backslash \{e_0\}$ and $0$ else. But, by definition of good spanning trees, it is immediate that:
$$\text{$T$ is a good spanning tree of $M$} \Longleftrightarrow \! \begin{array}{c} \text{$T$ is a good spanning tree of $M \backslash \{e_0\}$}\\  \text{and the first dart of $e_0$ appears} \\ \text{before its second in the tour of $T$.} \end{array}$$
But as $e_0 \in HE(C(e_0))$, the second condition of the right hand side is true by lemma \ref{lem_chboucle}. Finally, the coefficient of $T$ is $+1$ if $T$ is a good spanning subtree of $M$ and $0$ else.
\end{description}
\end{proof}

We are now ready to state our second main result: for this , we have to give a rooted map structure to our $G$. This is possible in multiple ways (choice of the map structure and of the place of the root).\\

\begin{Th}\label{th:N_combi}
The polynomial $N$ associated to the underlying graph $G$ of a rooted map $M$ is given by the following combinatorial formula:
\begin{equation}\label{eq:N_combi}
N(G)=\sum_{\substack{T \ \text{good spanning} \\ \text{ tree of }M}} \left[ \prod_{\substack{e \in HE(G) \\ e \notin T}} \big( x_{\alpha(e)}-x_{\omega(e)} \big) \right].
\end{equation}
\end{Th}

\begin{proof}
This is an immediate consequence of paragraph \ref{cons_rat_func} and Proposition \ref{prop:coef_arbres}.
\end{proof}

Of course, the good spanning trees depend on the map structure chosen on the graph $G$. So the theorem implies that the right member does not depend on it, which is quite surprising.\\

\section{A condition of factorization}\label{sectchainfact}

\subsection{Chain factorization}

In the previous section, we have given an additive formula for the numerator of the reduced fraction $\Psi_P$. Greene's formula for planar posets (see subsection \ref{subsect:background}) and the example of Figure \ref{fig:ex_factorization} show that, in some cases, it can also be written as a product of non-trivial factors. In this paragraph, we give a simple graphical condition on a graph $G$, which implies the factorization of $N(G)$.\\
Then, in the next paragraph we prove that, although our condition is not a necessary condition (see open problem \ref{subsubsect:openpb_facto}), it explains the fact that $N$ is a product of degree $1$ terms for planar posets.

 In this section, we will assume that all the graphs are connected, have no circuits and no transitivity relation (an edge going from the beginning to the end of a chain). As the value of $N$ on disconnected graphs is $0$ and Hasse diagrams of posets always fulfill the two others assumptions, we do not lose in generality. This means that, if we consider a chain $c$, there is no edges between the vertices of the chain except of course the edges of the chain itself.

\begin{figure}[h]
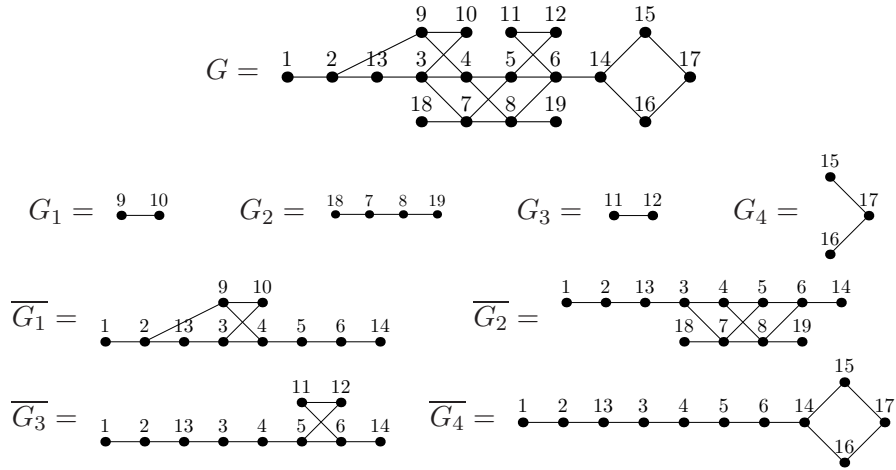

$$G=
\begin{array}{c}
\mbox{\epsfig{file=figure-56.ps,width=159pt}}
\end{array}
$$
$$G_1 =
\begin{array}{c}
\mbox{\epsfig{file=figure-57.ps,width=20pt}}
\end{array}
\qquad
G_2 =
\begin{array}{c}
\mbox{\epsfig{file=figure-58.ps,width=45pt}}
\end{array}
\qquad
G_3 =
\begin{array}{c}
\mbox{\epsfig{file=figure-59.ps,width=22pt}}
\end{array}
\qquad
G_4 =
\begin{array}{c}
\mbox{\epsfig{file=figure-60.ps,width=22pt}}
\end{array}$$
$$
\begin{array}{cc}

\overline{G_1}=
\begin{array}{c}
\mbox{\epsfig{file=figure-61.ps,width=110pt}}
\end{array}
&

\overline{G_2}=
\begin{array}{c}
\mbox{\epsfig{file=figure-62.ps,width=110pt}}
\end{array}
\\
\overline{G_3}=
\begin{array}{c}
\mbox{\epsfig{file=figure-63.ps,width=110pt}}
\end{array}
&
\overline{G_4}=
\begin{array}{c}
\mbox{\epsfig{file=figure-64.ps,width=143pt}}
\end{array}
\end{array}
$$
\caption{A graph $G$ with a chain $c$, the components $G_i$ of $G \setminus c$ and the corresponding regions $\overline{G_i}$.\label{fig:graphe_chaine_regions}} 
\end{figure}

Let $G$ be a graph, $c$ a chain of $G$, $V_c$ the set of vertices of $c$ (including the origin and te end of the chain) and $G_1$, \ldots, $G_k$ be all the connected component of $G \setminus V_c$. The complete subgraphs $\overline{G_i}=G_i \cup V_c$ (for $1\leq i \leq k$) will be called region of $G$.
Consider, for example, the graph of Figure \ref{fig:graphe_chaine_regions} and the chain $c=(1,2,13,3,4,5,6,14)$. In this case, the graph $G \setminus V_c$ has four connected components.

We can now state our third main result:
\begin{Th}
\label{factorization_chaine}
Let $G$ be a graph, $c$ a chain of $G$ and $\overline{G_1}, \overline{G_2}, \dots, \overline{G_k}$ be the corresponding regions of $G$.
Then one has:
$$
N(G) = \prod_{j=1}^{k} N(\overline{G_j} ).
$$
\end{Th}

For example, the numerator of the rational function associated to the graph of Figure \ref{fig:graphe_chaine_regions} can be written as a product of four non-trivial factors.\\

\emph{Proof.} The central idea is to apply Theorem \ref{th:boucle} page \pageref{th:boucle} on cycles $C$ contained in one region and such that $HE(C) \cap c = \emptyset$. This means that the edges of $c$ can appear in $C$, but only in the \emph{wrong} direction: so, when we apply Proposition \ref{th:boucle_N}, we do not cut the chain $c$.\\

\begin{figure}[h]
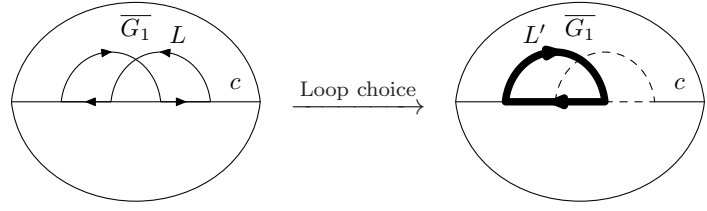

\begin{center}
$
\begin{CD}
\begin{array}{c} \mbox{\epsfig{file=figure-65.ps,width=96 pt}} \end{array}  @>\text{Loop choice}>> \begin{array}{c} \mbox{\epsfig{file=figure-66.ps,width=96 pt}} \end{array}\\
\end{CD}
$
\end{center}
\caption{\emph{Good} choice of cycle\label{fig:cycle_choice}}
\end{figure}

The first step is to prove the existence of such cycles. This is done in Lemma \ref{lemma:choose_cycle} (see Figure \ref{fig:cycle_choice} for an illustration).

\begin{Lem}
\label{lemma:choose_cycle}
Let $G$ be a graph and $c$ a chain of $G$. Denote by $\overline{G_1},\ldots,\overline{G_k}$ the corresponding regions.
If $\overline{G_1}$ is not a tree, there exists a cycle $C$ in $\overline{G_1}$ such that $HE(C) \cap c = \emptyset$.
\end{Lem}

\begin{proof}
Choose any cycle $C_0$ of $\overline{G_1}$. Two cases have to be examined:
\begin{enumerate}[1)]
\item The cycle $C_0$ has no vertices in common with $c$. Nothing has to be done.
\item The cycle $C_0=(e_1,\ldots,e_l)$ has at least one vertex in common with $c$. As a cycle is not transformed if one makes a cyclic permutation of its edges, one can assume that $V'=\extr_1(e_1)$ is a vertex of $c$. Let us denote by $h$ the smallest index such that $V''=\extr_2(e_h)$ is also a vertex of $c$ (it necessarily exists because $\extr_2(e_l)=\extr_1(e_1)$ is a vertex of $c$).
But there is a subchain (eventually empty) $c'$ of $c$ going from $V'$ to $V''$ (resp. from $V''$ to $V'$) if $V' \leq V''$ (resp. if $V'' \leq V'$). Now, we just have to define $C$ as:
$$
L= \left\{
\begin{array}{ll}
	(e_1,\ldots,e_l) \cdot \overline{c'} & \text{if }V' \leq V''\\
	(e_l,\ldots,e_1) \cdot \overline{c'}& \text{if } V''\le V'\\
\end{array} \right.,
$$
 where $\overline{c'}$ denotes the chain $c'$ in the other direction (this implies that all the edges of $c'$ are in the \emph{wrong} direction in $C$, so $HE(C) \cap c =\emptyset$). 
\end{enumerate}
\end{proof}

Let us come back to the proof of Theorem \ref{factorization_chaine}. We make a proof by induction on $k$. If $k=1$, then the result is trivial.\\

Suppose now that our proposition is true for $k=n-1$. Let $G$ be a graph and $c$ a chain of $G$, such that there are $n$ associated regions $\overline{G_1}, \dots, \overline{G_n}$.\\

If $\overline{G_1}$ is a tree, one can prune it 
to obtain the chain $c$. We can remove the same vertices and edges from the whole graph $G$ because the removed vertices are not linked with an other $G_i$ (as the $G_i$ are different connected components of $G \setminus V_c$). Thanks to the pruning-invariance lemma \ref{lem:emond} page \pageref{lem:emond}, one has:
$$N(G) = N(\bigcup_{i=2}^n \overline{G_i}) = \prod_{i=2}^k N(\overline{G_i}),$$
where the second equality is due to the induction hypothesis. The theorem is proved in this case.\\

If $\overline{G_1}$ is not a tree, we proceed by induction over the cyclomatic number of $\overline{G_1}$. 

Lemma \ref{lemma:choose_cycle} gives us a cycle $C_1$ of $\overline{G_1}$ such that $HE(C_1) \cap V_c = \emptyset$. Applying Proposition \ref{th:boucle_N} on $C$, one has
$$
N(G) = \sum_{{E_1 \subset HE(C_1) \atop E_1 \neq \emptyset}} \pm N(G \setminus E_1) \left(\prod_{e \in E_1}\left( x_{\alpha(e)} - x_{\omega(e)} \right) \right).
$$

Some of the graphs $G \setminus E_1$ are disconnected (if and only if $\overline{G_1} \setminus E_1$ is disconnected). The value of $N$ on these graphs is $0$. So they do not appear in the formulas \ref{preuve_fact_1} and \ref{preuve_fact_2}.\\

Each connected graph $G \setminus E_1$ contains the chain $c$ (thanks to the assumption $HE(C_1) \cap c = \emptyset$). The associated regions are $\overline{G_2},\ldots,\overline{G_n}$ and $\overline{G_1} \setminus E_1$ (the last region can in fact be a union of several regions but it does not matter). But $\overline{G_1} \setminus E_1$ has a strictly lower cyclomatic number than $\overline{G_1}$ so we can use the induction hypothesis
$$
N(G \setminus E_1) = N(\overline{G_1} \setminus E_1) \cdot N(\overline{G_2}) \cdot \ldots \cdot  N(\overline{G_n}).
$$
Finally,
\begin{equation}\label{preuve_fact_1}
N(G) = \left( \sum_{E_1 \subset HE(C_1) \atop E_1 \neq \emptyset} \pm  \prod_{e \in E_1} (x_{\alpha(e)}-x_{\omega(e)}) N(\overline{G_1} \setminus E_1) \right) \cdot N(\overline{G_2}) \cdot \ldots \cdot  N(\overline{G_n}),
\end{equation}
where the sum is restricted to the sets $E_1$ such that $G \setminus E_1$ is connected.
But we can use Proposition \ref{th:boucle_N} page \pageref{th:boucle_N}  with the same cycle $C$ in $G_1$:
\begin{equation}\label{preuve_fact_2}
N(\overline{G_1}) = \sum_{E_1 \subset HE(C_1) \atop E_1 \neq \emptyset } \pm  \prod_{e \in E_1} (x_{\alpha(e)}-x_{\omega(e)}) N(\overline{G_1} \setminus E_1) ,
\end{equation}
where the sum is also restricted to the sets $E_1$ such that $\overline{G_1} \setminus E_1$ is connected, or equivalently such that $G \setminus E_1$ is connected.\\

This ends the proof of Theorem \ref{factorization_chaine}.

\subsection{Complete factorization of planar posets}\label{green_theorem}

In his paper \cite{GREENE}, C. Greene has given a closed expression for the sum $\Psi(G)$ when $G$ is the minimal graph (Hasse diagram) of a \emph{planar} poset (Theorem \ref{th:greene}). In this case, the numerator $N(G)$ can be written as a product of terms of degree $1$ (Theorem \ref{th:greene}). We will see that this factorization property is a consequence of Theorem \ref{factorization_chaine} and give a new proof of Greene's Theorem.\\

Let us begin by defining precisely planar posets:

\begin{Def} \label{def:planar}
We will say that the drawing of an oriented graph (without circuit) is ordered-embedded in $\mathbb{R} \times \mathbb{R}$ if
\begin{itemize}
 \item the origin of an edge is always at the left of its end ;
 \item the edges are straight lines.
\end{itemize}
A graph $G$ is said planar if it can be ordered embedded  in $\mathbb{R} \times \mathbb{R}$ without edge-crossings.
If $G$ is a graph, we denote by $G_{0,\infty}$ the graph obtained from $G$ by adding:
\begin{itemize}
 \item A vertex $0$ (called \emph{minimal} vertex) and, for each vertex $v$ of $G$ which is not the end of any edge of $G$, an edge going from $0$ to $v$.
 \item A vertex $\infty$ (called \emph{maximal} vertex) and, for each vertex $v$ of $G$ which is not the origin of any edge of $G$, an edge going from $v$ to $\infty$.
\end{itemize}
A graph $G$ is said strongly planar if the graph $G_{0,\infty}$ is planar.\\
A poset $\mathcal{P}$ is \emph{planar} if its minimal graph $G$ is strongly planar.
\end{Def}

\begin{figure}[h]
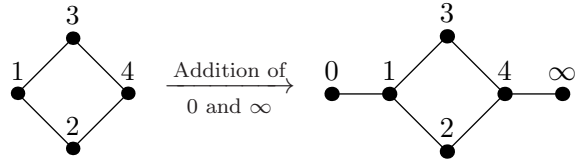

\begin{center}
$
\begin{CD}
\begin{array}{c} \mbox{\epsfig{file=figure-10.ps,width=48 pt}} \end{array}  @>\text{Addition of}>0\text{ and }\infty> \begin{array}{c} \mbox{\epsfig{file=figure-67.ps,width=96 pt}} \end{array}\\
\end{CD}
$
\end{center}
\caption{The graph $G$ is strongly planar.\label{fig:planar_poset}}
\end{figure}

\begin{figure}[h]
\begin{center}
$
\begin{CD}
\begin{array}{c} \mbox{\epsfig{file=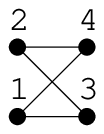,width=24 pt}} \end{array}  @= \begin{array}{c} \mbox{\epsfig{file=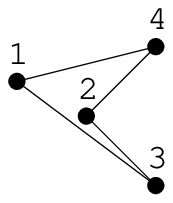,width=48 pt}} \end{array}\\
@V\text{Addition of }0\text{ and }\infty VV  @VV\text{Addition of }0\text{ and }\infty V\\
\begin{array}{c} \mbox{\epsfig{file=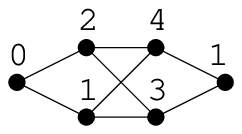,width=72 pt}} \end{array} @= \begin{array}{c} \mbox{\epsfig{file=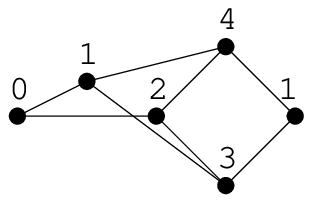,width=96 pt}} \end{array}
\end{CD}
$
\end{center}
\caption{The graph $G$ is planar, but not strongly planar.\label{fig:non_planar_poset}}
\end{figure}

Almost all drawings of this paper (except in section \ref{sect:combi}) are ordered-embedded in $\mathbb{R} \times \mathbb{R}$. See Figure  \ref{fig:planar_poset} and \ref{fig:non_planar_poset} for examples of strongly planar and non strongly planar graphs.

Note that the induced graph of a strongly planar graph is strongly planar (note that, however, if we erase some edges, we can obtain a non strongly planar graph). In particular, the regions of a strongly planar graph with respect to a chain are the graph of strongly planar graphs.\\

\begin{figure}[h]
\begin{center}
$
\mbox{\epsfig{file=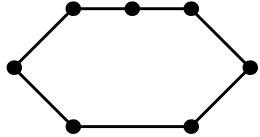,width=96 pt}}
$
\end{center}
\caption{A diamond \label{diamond}}
\end{figure}

 Moreover, a graph with one cycle and without vertices with valence $1$ is strongly planar if and only if it has a unique maximal and a unique minimal element. In this case, we will call it a diamond (an example is drawn on Figure \ref{diamond}).\\

These definitions are relevant because there is a closed formula for $\Psi_{\mathcal{P}}$ for planar posets:

\begin{Th}[Greene \cite{GREENE}]
\label{th:greene}
Let $P$ be a planar poset, then:
$$
\Psi_{\mathcal{P}} =
\left\{
\begin{array}{cl}
0 & \text{if } P \text{ is not connected;} \\
\prod_{y,z \in P} (x_y-x_z)^{\mu_P(y,z)} & \text{if } P \text{ is connected,}
\end{array}
\right. ,
$$
\noindent where $\mu_P(x,y)$ denotes the M\"{o}bius function of the poset $P$.
\end{Th}

We will show that we can find disconnecting chains in any strongly planar graphs, explaining the fact that the function $N(G)$ can be factorized into factors of degree $1$.\\

\begin{Prop}
\label{prop:fractionNonTrivial}
Let $G$ be a strongly planar oriented graph with a number of cycle greater than $1$, then there is a chain of $G$, separating $G$ in two non-trivial regions (each region contains at least one cycle).
\end{Prop}

\begin{proof}
Eventually by pruning it, one can assume that $G$ has no vertices with valence $1$. As it has at least two cycles, it has one vertex $c_2$ of valence $3$ or more. So, up to a left-right symmetry, we are in one of the two following cases (in the second case, we assume that $c_2$ is the end of \emph{exactly} 2 edges).

\begin{multicols}{2}
\begin{center}
\epsfig{file=figure-73.ps,width=26 pt}

\epsfig{file=figure-74.ps,width=48 pt}
\end{center}
\end{multicols}

In the first case, let us label the vertices as below:

\begin{center}
\epsfig{file=figure-75.ps,width=30 pt}
\end{center}

In the second case, we define by induction $c_i$ for $i\geq 3$: we choose for $c_i$ any vertex such that there is an edge of origin $c_{i-1}$ and of end $c_i$. For a $k \geq 3$, one can not define $c_{k+1}$ if $c_k$ is not the origin of any edge. Then, as $c_{k}$ is not a vertex of valence $1$, it is the end of an edge coming from a vertex $b \neq c_{k-1}$. Finally, we call $c_1$ and $a$ the origins of the two edges whose ends are $c_2$: which one is $c_1$ and which one is $a$ depends of whether $b$ is above or below $c_{k-1}$ (see the figure below).

\begin{multicols}{2}
\begin{center}
\epsfig{file=figure-76.ps,width=48 pt}

\epsfig{file=figure-77.ps,width=48 pt}
\end{center}
\end{multicols}

In every case, the $c_i$ are the vertices of a chain $c$ of $G$, which can be extended to a maximal chain $c_{\text{max}}$. Recall that with Greene's definition of a planar graph, the graph $G_{0,\infty}$, \text{i.e.} can still be ordered-embedded in the plan. Then there is a chain in $G_{0,\infty}$ containing $v_0$, $c_{\text{max}}$ and $v_\infty$. It splits $G_{0,\infty}$ into at least two regions, one containing $a$ and one containing $b$. The same is true for the chain $c_{\text{max}}$ in $G$. But, as $G$ has no vertices, the corresponding regions have at least one cycle.
\end{proof}

\begin{Corol}
Let $G$ be a connected strongly planar poset. By iterating chain factorization, one can write $N(G)$ as a product of numerators of rational functions associated to diamonds.
 \label{corol:factorization:diamond}
\end{Corol}

\begin{proof}
 Proposition \ref{prop:fractionNonTrivial} page \pageref{prop:fractionNonTrivial} and Theorem \ref{factorization_chaine} page \pageref{factorization_chaine} imply that $N(G)$ can be factorized as the product of numerators of subgraphs with one cycle. As these subgraphs are strongly planar, after pruning, they are diamonds, which ends the proof.
\end{proof}

Note that for a diamond, the function $N$ has a closed expression (paragraph \ref{subsect:ex_1cycle}):
$$N(D)=x_{min(D)} - x_{max(D)},$$
or equivalently,
$$\Psi(D)= \prod_{y,z \in P} (x_y-x_z)^{\mu_D(y,z)},$$
where $\mu_D$ is the Möbius function of the poset associated to the diamond $D$.\\

The last property can be extended to any planar poset thanks to the following compatibility between disconnecting chain and Möbius function:

\begin{Prop}
\label{prop:factorization:mobius}
Let $P$ be a poset, $c$ a chain of the Hasse diagram of $P$ (\textit{i.e.} the minimal graph representing $P$), $P_1, \dots, P_n$ the $n$ region associated with $c$, and $i,j$ two different elements  of $P$, then
$$
\mu_P(i,j) = \left\{
\begin{array}{ll}
-1 & \text{if }i \preceq j,\\
\sum_{k=1}^{n} \mu_{P_k}(i,j) & \text{otherwise}.
\end{array}
\right.
$$
We assume that $\mu_Q(i,j)=0$ if $i \not\in Q$ or $j \not\in Q$.
\end{Prop}

The proof is postponed to paragraph \ref{subsect:proof_chain_mobius}.\\

This proposition together with corollary \ref{corol:factorization:diamond} proves Greene's theorem. In fact, this proof works also for some non-planar posets (and hence Greene's formula is true for these posets). For example, the poset of the figure \ref{facto_chaine} is not planar but can be factorised and the numerator can be expressed with the Möbius function: this is the case of any gluing of diamonds along chains.

\begin{figure}[h]
$$
\begin{array}{rl}
N\left(\begin{array}{c}
\mbox{\epsfig{file=figure-78.ps,width=88pt}}
\end{array}\right)
&
=
N \left( \begin{array}{c}  
\mbox{\epsfig{file=figure-79.ps,width=53pt}}
\end{array} \right)
\cdot
N \left( \begin{array}{c}
\mbox{\epsfig{file=figure-80.ps,width=53pt}}
\end{array} \right)\\
& \qquad
\cdot \
N \left( \begin{array}{c}
\mbox{\epsfig{file=figure-81.ps,width=53pt}}
\end{array} \right)
\\
&= (x_1-x_4) . (x_2-x_5) . (x_3-x_6)
\end{array}
$$
\caption{A non-planar poset for which Greene's formula is true.  \label{facto_chaine}}
\end{figure}

\subsection{Chain and Möbius function}\label{subsect:proof_chain_mobius}
This paragraph is the proof of the technical Proposition \ref{prop:factorization:mobius}

\begin{proof}

When $i \preceq j$ (there is an edge from $i$ to $j$ in the Hasse diagram of the poset), one always has $\mu_P(i,j) = -1$.\\

When $i \leq j$, but $i \not\preceq j$, four cases have to be examined:
\begin{description}
 \item[first case] $i,j$ do not belong to $V_c$ and in different regions of the poset;
 \item[second case] $i,j$ do not belong to $V_c$, but are in the same region of the poset ;
 \item[third case]  $i$ is an element of $V_c$, but $j$ is not ;
 \item[fourth case] $i$ and $j$ are two elements of $V_c$.
\end{description}

Figure \ref{fig:mobius:case:1}, \ref{fig:mobius:case:2}, \ref{fig:mobius:case:3} and  \ref{fig:mobius:case:4} summarize the four cases. Note that the case where $i$ does not belong to $V_c$, but $j$ does, can be obtained from the third one by considering the opposite poset.\\

Let $P_1, \dots, P_n$ be the $n$ regions associated with $P$.\\

We denote by $[a,b]_P$ the set
$$
[a,b]_P = \{ k | a \le_P k \le_P b \},
$$
and by $[a,b[_P$ the set
$$
[a,b[_P = \{ k | a \le_P k <_P b \}.$$

Note that $[i,j]_{P_1} = [i,j]_P \cap P_1$. This property is not true for any poset associated to a complete subgraph of $G$, the fact that $P_1$ is a region defined by a disconnecting chain is here very important.\\

If $[i,j]_P$ has a non-empty intersection with $V_c$, we denote by $L$ the maximal element of this intersection.

\comment{We denote by $l$ the integer defined by:
\begin{multline*}
 l = max \{ length(v) | v \text{ is a chain of } P \\
\text{such that the first element of } v \text{ are in } [i,L]_P \\
\text{ and the other are in }P_1\setminus[i,j]_P \}.
\end{multline*}

For example, in the figure \ref{fig:mobius:1}, one has $[i;j]_P \cap V_c = \{3,4,5\}$, so $L$ exists and is the vertex labeled $5$
.

\begin{figure}[h]
\begin{center}
$
\mbox{\epsfig{file=figure-82.ps,width=140 pt}}
$
\end{center}
\caption{Illustration of the definition of $L$.\label{fig:mobius:1}}
\end{figure}
}

\begin{enumerate}[1)]
\item \label{en:case:1} Suppose that $i \in P_2 \setminus V_c$ and $j \in P_1 \setminus V_c$. We want to prove that $\mu_P(i,j)=0$ and we assume (proof by induction) that it is true for any $j' \in P_1 \setminus V_c$ such that $j' < j$.

\begin{figure}[h]
\begin{center}
$
\mbox{\epsfig{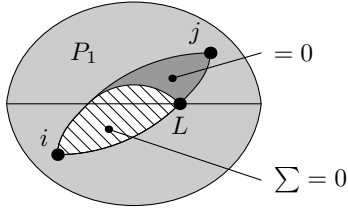}}
$
\end{center}
\caption{Case \ref{en:case:1}: $i \not\in V_c$ and $j$ is not in the same region than $j$. \label{fig:mobius:case:1}}
\end{figure}

As $i \leq j$, there is a chain in the Hasse diagram of $P$ going from $i$ to $j$. As $c$ is a chain separating $P_1$ and $P_2$, any chain from $i$ to $j$ intersect $V_c$. Thus $L$ exists and any element between $i$ and $j$ which is not in $P_1$, is lower or equal to $L$. So
$$
[i,j]_P \cap (P_2 \cup \dots \cup P_m) \subseteq [i,L]_P \subseteq [i,j]_P.
$$
By definition of the Möbius function we obtain,
$$
\mu_P(i,j) = -\sum_{k \in [i,L]_P} \mu_{P}(i,k) - \sum_{k \in [i,j[_{P} \cap P_1 \setminus  [i,L]_P} \mu_{P}(i,k)
$$
As 
$$
\sum_{k \in [i,L]_P} \mu_{P}(i,k) = 0
$$
one has:
\begin{equation}
\label{eq:mobius}
\mu_P(i,j) = - \sum_{k \in [i,j[_P \cap P_1 \setminus [i,L]_P} \mu_{P}(i,k)
\end{equation}
By induction hypothesis, $\mu_P(i,j)=0$.

\item \label{en:case:2} Suppose now that $i,j \in P_1 \setminus V_c$. We want to prove that $\mu_P(i,j)=\mu_{i,j}(P_1)$

\begin{figure}[h]
\begin{center}
$
\mbox{\epsfig{file=figure-84.ps,width=130 pt}}
$
\end{center}
\caption{Case \ref{en:case:2}: $i \not\in V_c$ and $j$ is in the same region than $i$. \label{fig:mobius:case:2}}
\end{figure}

By definition of the Möbius function, we have,
$$
\mu_P(i,j) = - \sum_{k \in [i,j[_P \cap (P_2 \cup \dots \cup P_m \setminus V_c)}\mu_P(i,j) - \sum_{k \in [i,j[_{P_1}} \mu_P(i,j).
$$
The case $1$ gives: $\sum_{k \in [i,j[_P \cap (P_2 \cup \dots \cup P_m \setminus V_c)}\mu_P(i,j) =0$. Therefore,
\begin{equation}
\label{eq:mobius:2}
\mu_P(i,j) = \sum_{k\in[i,j[_{P_1}} \mu_P(i,j).
\end{equation}
and an immediate induction proves that $\mu_P(i,j)=\mu_{i,j}(P_1)$.

\item \label{en:case:3} Suppose that $i \in V_c$ and $j \in P_1 \setminus V_c$. As $i \in V_c \cap [i,j]_P$, the set is not empty and $L$ exists.

\begin{figure}[h]
\begin{center}
$
\mbox{\epsfig{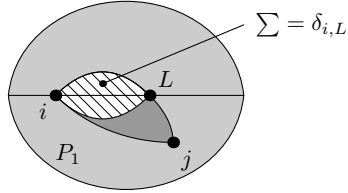}}
$
\end{center}
\caption{Case \ref{en:case:3}: $i \in V_c$ and $j \not\in V_c$. \label{fig:mobius:case:3}}
\end{figure}

We will prove now that $\mu_P(i,j) = \mu_{i,j}(P_1)$ by induction on $j$.
As
$$\sum_{k \in [i,L]_P} \mu_{P}(i,k)=\delta_{i,L}$$
one has:
$$
\mu_P(i,j) = - \sum_{k \in [i,j[_P \setminus [i,L]_P} \mu_{P}(i,k)
 - \delta_{i,L}
$$
Similarly,
$$\mu_{i,j}(P_1)= - \sum_{k \in [i,j[_{P_1} \setminus [i,L]_P} \mu_{i,k}(P_1) - \delta_{i,L}$$
But $[i,j[_P \setminus [i,L]_P = [i,j[_{P_1} \setminus [i,L]_P$ (see the proof of case $1$), so an immediate induction on $j$ finishes the proof in this case.

\comment{
\begin{enumerate}[i)]
\item If $l=1$, two cases have to be examined. If $i \preceq j$ then $\mu_P(i,j) = \mu_{i,j}(P_1)=-1$.

If $i \not\preceq j$ then $L \preceq i$ and, in the Hasse diagram of the poset induced by $[i,j]_P$, $i$ is connected on the poset by $L$ (see figure \ref{fig:mobius:case:3:l1}).

\begin{figure}[h]
\begin{center}
$
\mbox{\epsfig{file=figure-86.ps,width=110 pt}}
$
\end{center}
\caption{The subposet induced by $[i,j]_P$ when $l=1$ and $L \preceq i$  \label{fig:mobius:case:3:l1}}
\end{figure}

We deduce that $\mu_P(i,j) = \mu_{i,j}(P_1)=0$.
\item Suppose that $\mu_P(i,j) = \mu_{i,j}(P_1)$  for $l \le n$. Let $l=n+1$ then
$$
\mu_P(i,j)= - \sum_{k\in[i,j[_{P_1}} \mu_{P}(i,k) = - \sum_{k\in[i,j[_{P_1}} \mu_{i,k}(P_1).
$$
We deduce that $\mu_P(i,j) = \mu_{i,j}(P_1)$.
\end{enumerate}
}

\item \label{en:case:4} Suppose that $i \in V_c$ and $j \in V_c$. We want to prove by induction on $j$ ($i \not\preceq j$) that $\mu_P(i,j)=\sum\limits_{l=1}^n \mu_{P_l}(i,j)$.

\begin{figure}[h]
\begin{center}
$
\mbox{\epsfig{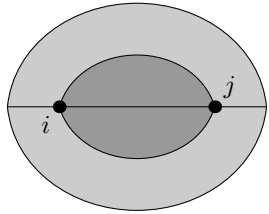}}
$
\end{center}
\caption{Case \ref{en:case:4}: $i$ and $j$ belong to $V_c$.\label{fig:mobius:case:4}}
\end{figure}

By definition of the Möbius function, we have
$$
\mu_P(i,j) 	= - \sum_{k \in [i,j[_P} \mu_{P}(i,k)
$$
Using case \ref{en:case:3} of this poof and the induction hypothesis, we know that $\mu_{P}(i,k)=\sum_{l=1}^n \mu_{P_l}(i,k)$ if $k \in [i,j[_P$ except for:
\begin{description}
\item[$\mathbf{k=i}$] In this case, $\mu_{P}(i,i)=\mu_{P_l}(i,i)=1$, thus $\mu_{P}(i,i)=\sum\limits_{l=1}^n \mu_{P_l}(i,i) - (n-1)$.
\item[$\mathbf{k=i_1}$] where $i_1$ is defined by $i_1 \in V_c$ and $i \preceq i_1$. In this case, one has $\mu_{P}(i,i_1)=\mu_{P_l}(i,i_1)=-1$, thus $\mu_{P}(i,i_1)=\sum\limits_{l=1}^n \mu_{P_l}(i,i_1) + (n-1)$.
\end{description}
Finally, one has:
$$
\mu_P(i,j) = \sum_{l=1}^n\left( - \sum_{k \in [i,j[_P \setminus c} \mu_{P_l}(i,k) \right) - (n-1) + (n-1).
$$
Using the definition of the Möbius function for the $P_l$, this ends the proof of the proposition.
\end{enumerate}

\end{proof}

\section*{Acknowledgements}
\label{sec:ack}
The authors are grateful to A. Lascoux for his suggestion to work on this rational functions.

\nocite{*}
\bibliographystyle{alpha}
\bibliography{factorisation}

\end{document}